\documentclass[12pt,reqno]{amsart}

\usepackage{amsmath,amsthm,amssymb,mathtools}
\usepackage{stmaryrd,xcolor}
\SetSymbolFont{stmry}{bold}{U}{stmry}{m}{n}
\usepackage{enumitem}
\usepackage{hyperref}
\usepackage{cleveref}
\usepackage{geometry}
\usepackage[backend=biber,style=numeric,sorting=nyt,giveninits=true]{biblatex}
\crefformat{section}{\S#2#1#3}
\crefformat{subsection}{\S#2#1#3}
\crefformat{subsubsection}{\S#2#1#3}

\allowdisplaybreaks

\newtheorem{theorem}{Theorem}[section]
\newtheorem{lemma}[theorem]{Lemma}
\newtheorem{corollary}[theorem]{Corollary}

\theoremstyle{definition}

\theoremstyle{remark}

\numberwithin{equation}{section}

\newcommand{\floor}[1]{\left \lfloor #1 \right \rfloor}

\def\bfh{{\mathbf h}} 

\def\bfj{{\mathbf j}}
\def\bfk{{\mathbf k}}
 
\def\bfm{{\mathbf m}} 
\def\bfn{{\mathbf n}}

\def\bfr{{\mathbf r}}
\def\bfs{{\mathbf s}}
\def\bft{{\mathbf t}}

\def\bfv{{\mathbf v}}

\def\bfx{{\mathbf x}}

\def\calB{{\mathcal B}}

\def\calD{{\mathcal D}}

\def\C{{\mathbb C}}\def\D{{\mathbb D}}

\def\N{{\mathbb N}}  
\def\R{{\mathbb R}}
\def\Z{{\mathbb Z}}\def\Q{{\mathbb Q}}

\def\a{{\alpha}}

\def\eps{\varepsilon}

\def\le{\leqslant} \def\ge{\geqslant}

\def\d{{\,{\rm d}}}

\theoremstyle{definition}

\theoremstyle{remark}
\crefname{theorem}{Theorem}{Theorems}
\crefname{lemma}{Lemma}{Lemmas}
\crefname{proposition}{Proposition}{Propositions}
\crefname{corollary}{Corollary}{Corollaries}
\crefname{definition}{Definition}{Definitions}
\crefname{remark}{Remark}{Remarks}
\crefname{conjecture}{Conjecture}{Conjectures}
\Crefname{conjecture}{Conjecture}{Conjectures}
\crefname{problem}{Problem}{Problems}
\Crefname{problem}{Problem}{Problems}
\crefname{hypothesis}{Hypothesis}{Hypotheses}
\Crefname{hypothesis}{Hypothesis}{Hypotheses}

\newcommand{\e}{\mathrm e}

\newcommand{\bfmu}{\boldsymbol\mu}
\newcommand{\bflam}{\boldsymbol\lambda}

\newcommand{\bfzero}{\mathbf 0}
\newcommand{\dd}{\,\mathrm d}
\newcommand{\disc}{\operatorname{disc}}

\begin{document}
\title[Bounds for the binary quartic Weyl sum]{Bounds for the binary quartic Weyl sum}
\author[Daniel Flores Galiote]{Daniel Flores Galiote}
\address{Department of Mathematics, Rice University}
\email{df66@rice.edu}
\date{}
\hypersetup{
  pdftitle={Bounds for the binary quartic Weyl sum},
  pdfauthor={Daniel Flores Galiote}
}

\begin{abstract} 
We prove that
\[
\sum_{1\le x,y\le N}\e(\a F(x,y))
\ll_{F,\a,\eps}N^{5/3+\eps}
\]
for every $\eps>0$, real quadratic irrational $\a$,
and integral binary quartic form $F$ that is not the
fourth power of a linear form over $\C$.
The proof uses the two-dimensional $q$-analogue of van der
Corput's method, following Heath-Brown's work in the one-variable case.
The heart of our argument rests on new bounds for complete
binary cubic Weyl sums with linear terms, which retain
the dependence on the cubic discriminant and the linear
coefficients.
\end{abstract}

\maketitle

% The order follows suggestedorder.tex; the cubic-sum proofs come last.

\section{Introduction}\label{sec:setup}

For over a century now, the Hardy-Littlewood circle method has 
been a driving force behind many amazing results in analytic 
number theory. Its applications include Waring's problem
\cite{Hardy1920}, the ternary Goldbach problem
\cite{Vinogradov1937,Helfgott2013}, and Birch's Theorem on
forms in many variables \cite{Birch1962}. More recently,
Maynard \cite{Maynard2019} used the circle method together
with sieve methods to prove that there are infinitely many
primes omitting any one prescribed decimal digit.

When the circle method is applied to questions concerning the
solubility of polynomial equations over the integers, its success
often hinges on establishing adequate bounds for Weyl sums.
A basic example is the one-variable Weyl sum
\[
S_d(\a;N)=\sum_{1\le n\le N}\e(\a n^d),\qquad \e(z):=\exp(2\pi iz),
\]
where $d\ge2$ is an integer, $\a\in\R$, and $N\ge2$ is an integer.

The method of Weyl \cite{Weyl1916}, developed quantitatively
by Hardy and Littlewood \cite{Hardy1920} and later by 
Vinogradov \cite{Vinogradov1947}, repeatedly applies
Cauchy's inequality along with a differencing procedure to
reduce the degree of the polynomial until it is linear in $n$.
This leads to the bound
\[S_d(\a;N) \ll_d N^{1+\eps}(q^{-1}+N^{-1}+qN^{-d})^{2^{1-d}},\]
whenever there exist coprime integers $a\in\Z$ and $q\in\N$ such that
$|\a - a/q| < 1/q^2$.

Thus, when $\a \in \R \backslash \Q$ is an algebraic number one 
may deduce from  Roth's Theorem \cite{Roth1955} the bound
\[
S_d(\a;N)\ll_{d,\a}N^{1-2^{1-d}+\eps}.
\]
For $d=2$ this gives essentially optimal square-root cancellation.
The exponents $3/4$ and $7/8$ in degrees three and four
remain the best available for arbitrary fixed algebraic
irrational arguments. The sharp mean-value estimates of Wooley
in degree three \cite{Wooley2016} and of Bourgain, Demeter and
Guth in higher degrees \cite{Bourgain2016} resolve the main
conjecture in Vinogradov's mean value theorem. In higher degrees,
these estimates also yield stronger bounds for Weyl sums.

For sums involving general forms in several variables, bounds
of comparable strength are much less widely available.
Integral binary quadratic forms with nonzero discriminant provide
a simple case. Completing the square reduces the required estimates 
to those for quadratic Weyl sums in one variable, giving the 
exponent $1+\eps$ for a sum over a box $[1,N]^2 \cap \Z^2$ when 
$\a$ is a fixed real algebraic irrational number.

In degree three, Heath-Brown \cite{Heath2010} improves the
one-variable exponent to $5/7$ for quadratic irrational
arguments under the assumption of the $abc$-conjecture. For binary
cubic forms with nonzero discriminant, Chowla and Davenport
\cite{Chowla1961} established a Weyl inequality giving the
exponent $3/2$ at algebraic irrational arguments.
Estimates for binary cubic sums also underlie Br\"udern and
Wooley's work \cite{Bruedern1998} on representations by sums
of binary cubic forms. Wooley \cite{Wooley1999,Wooley2000}
subsequently obtained the same exponent for every integral
binary cubic form that is not the cube of a linear form
over $\C$.

For quartic sums, Heath-Brown \cite[Theorem~1]{Heath2024}
proves unconditionally that
\[
S_4(\a;N)\ll_{\a,\eps}N^{5/6+\eps}
\]
when $\a$ is a real quadratic irrational. His proof uses the
$q$-analogue of van der Corput's method and rational approximations
whose denominators have factors of suitable sizes.

This paper concerns estimates for the binary quartic Weyl sum
\[
S_F(\a;N)=\sum_{1\le x,y\le N}\e(\a F(x,y)),
\]
where $F\in\Z[x,y]$ is a binary quartic form, $\a\in\R$,
and $N\ge2$.
Following Wooley \cite{Wooley1999,Wooley2000}, we call $F$
non-degenerate when it is not the fourth power of a linear form
over $\C$.

For non-degenerate $F$, Wooley's binary version of Weyl's
inequality \cite[Theorem~1]{Wooley1999}, combined with Roth's
Theorem \cite{Roth1955}, gives
\[
S_F(\a;N)\ll_{F,\a,\eps}N^{7/4+\eps}
\]
whenever $\a$ is a real algebraic irrational.
We improve this estimate when $\a$ is a quadratic irrational.

\begin{theorem}\label{thm:main}
For every non-degenerate binary quartic form $F\in\Z[x,y]$,
every real quadratic irrational $\a$,
and every $\eps>0$,
\[
S_F(\a;N)\ll_{F,\a,\eps}N^{5/3+\eps}.
\]
\end{theorem}

For the diagonal form $F(x,y)=x^4+y^4$, Heath-Brown's bound gives
\[
S_{x^4+y^4}(\a;N)
=\left(\sum_{1\le x\le N}\e(\a x^4)\right)^2
\ll_{\a,\eps}N^{5/3+\eps}.
\]
Thus \cref{thm:main} extends the exponent supplied by his
result for this diagonal form to every non-degenerate binary
quartic.

Our assumption permits repeated linear factors and does not
require $\disc(F)\ne0$. For the excluded form $F=x^4$, one has
\[
S_F(\a;N)=N\sum_{1\le x\le N}\e(\a x^4),
\]
so our bound would require the stronger one-variable exponent $2/3$.

\subsection{Overview of the proof}\label{subsec:overview}

Our proof uses the two-dimensional $q$-analogue of van der
Corput's method, following Heath-Brown's one-variable quartic argument
\cite{Heath2024}.
We first approximate $\a$ by $a/q$, with
\[q\asymp N^2, \quad \text{and} \quad |\a-a/q|\ll N^{-4},\]
which is possible for any badly approximable irrational.
This essentially reduces the problem to the case where $\a =a/q$.
We also require a factorization $q=q_1q_2$ in which $q_1$
has a prescribed size and $q_2$ avoids a fixed finite set of
primes. This is where we use the assumption that $\a$ is a 
quadratic irrational: the approximation results in
\cite{Heath2010,Heath2024} provide such denominators;
see \cref{sec:preliminaries}.

The $A$-process in \cref{sec:processes} differences our quartic Weyl sum
using shifts which are multiples of $q_1$. This leads one to
the problem of establishing bounds for cubic
Weyl sums modulo $q_2$, whose length relative to the modulus
can be optimized through the choice of $q_1$. The leading
terms of the polynomial arguments are then a constant multiple of
\[
G_{\bfh}(\bfm)=h_1F_x(\bfm)+h_2F_y(\bfm),
\qquad \bfh=(h_1,h_2)\in\Z^2,
\]
with the remaining terms linear in $\bfm$.
The $B$-process in \cref{sec:completion} expresses each
cubic Weyl sum as a weighted sum of complete cubic Weyl sums with
additional linear terms. To control this weighted average, we
require bounds for the individual complete Weyl sums that
retain the dependence on the linear coefficients. These bounds
are provided by our main technical result, \cref{thm:complete}.

For an integral binary cubic form $G$, an integer $r\ge1$,
and $\bfmu\in\Z^2$, write $\e_r(z)=\e(z/r)$ and put
\begin{equation}\label{eq:complete-def}
\Sigma(r;G,\bfmu)=
\sum_{\bfs\in(\Z/r\Z)^2}\e_r(G(\bfs)+\bfmu\cdot\bfs).
\end{equation}
For each fixed cubic $G$, we must bound a weighted sum of
$|\Sigma(r;G,\bfmu)|$ as $\bfmu$ ranges over translated rectangles.
Browning and Heath-Brown \cite{Browning2009} also estimate
quartic sums by differencing once and applying Poisson
summation to the resulting cubic sums. Their estimates retain
the dependence on the cubic coefficients and are averaged
over both the added linear coefficients and the shifts.
Here, our bounds retain the dependence on the linear coefficients. 
This allows us to count the choices for which the bound is large 
and control their total contribution after completion.

For one-variable sums, Cochrane and Zheng
\cite[Theorem~1.1]{Cochrane2000} treat the polynomial
$Ax^d+vx$ in arbitrary degree $d\ge2$. In particular, for
integers $A,v$, an odd prime $p$ with $p\nmid Ad(d-1)$,
and $k\ge2$, they prove
\[
\left|\sum_{x\bmod p^k}\e_{p^k}(Ax^d+vx)\right|
\le(d-1,p-1)p^{k/2}(p^k,v)^{1/2}.
\]
Their full Theorem also treats primes dividing $A$ or
$d(d-1)$, with additional factors involving their valuations.
For binary cubics, we require bounds that retain the dependence
on both added linear coefficients.

The technical heart of our argument is \cref{thm:complete}.
At prime-power moduli, our bounds for \eqref{eq:complete-def}
depend on the number of solutions of the first derivative congruences
\[
G_x(\bfs)+\mu_1\equiv G_y(\bfs)+\mu_2\equiv0\pmod{p^j}.
\]
When the Hessian determinant vanishes at a solution modulo $p$,
there may be many solutions at higher powers of $p$, allowing
larger values of the Weyl sum. The case $\bfmu=0$ illustrates
this difficulty: the first and second derivatives of $G$ all
vanish at the origin, and, for example,
\[
\Sigma(p^3;x^3+y^3,0)=p^4\qquad(p\ge5).
\]
We retain the dependence on the linear coefficients, so that
we can count the choices responsible for the larger bounds.
We must also retain the dependence on the cubic discriminant,
since the cubic forms arising from the $A$-process, $G_{\bfh}(\bfm)$, depend on
the shifts $\bfh$.

We identify the linear coefficients for which the Hessian
determinant can vanish at a solution of the derivative congruences
by a quadratic form. For
$G(x,y)=ax^3+bx^2y+cxy^2+dy^3$, write $\Delta=\disc(G)$
and put
\begin{equation}\label{eq:DG}
Q_G(\mu_1,\mu_2)=(c^2-3bd)\mu_1^2-(bc-9ad)\mu_1\mu_2
 +(b^2-3ac)\mu_2^2.
\end{equation}
The identity
\[
16Q_G(G_x,G_y)=(G_{xx}G_{yy}-G_{xy}^2)^2
\]
shows that, for $p\ge5$, a solution of
$G_x+\mu_1\equiv G_y+\mu_2\equiv0\pmod p$ at which the
Hessian determinant vanishes must satisfy
$Q_G(\bfmu)\equiv0\pmod p$. This gives a congruence in the
linear coefficients that we can count.
When $\Delta\ne0$, define
\[
c_G=\gcd(b^2-3ac,\,bc-9ad,\,c^2-3bd)>0,\qquad D_G=Q_G/c_G.
\]
We write $(r,\bfmu)=(r,\mu_1,\mu_2)$.

\begin{theorem}\label{thm:complete}
Let $r\ge1$ be an integer, let $\bfmu\in\Z^2$, and let
$G(x,y)=ax^3+bx^2y+cxy^2+dy^3$ be a nonzero integral binary
cubic form.
For every $\eps>0$, the following bounds hold.
\begin{enumerate}[label=\textup{(\roman*)},leftmargin=*]
\item If $\Delta\ne0$ and $(r,6)=1$, then
\begin{equation}\label{eq:complete-bound}
|\Sigma(r;G,\bfmu)|\ll
r^{1+\eps}(r^4,\Delta)^{1/4}
\min\left\{r^{1/4}(r,\bfmu)^{1/2},
(r,D_G(\bfmu))^{1/2}\right\}.
\end{equation}
\item If $\Delta=0$, one can choose a two by two integer matrix $V$
with $\det V\ne0$ and an integer $A\ne0$ such that $G\circ V$ is
either $AX^2Y$ or $AX^3$.  Write
\[
V^{\mathsf T}\bfmu
=\bigl(M_G(\bfmu),L_G(\bfmu)\bigr)^{\mathsf T}.
\]
For $(r,6A\det V)=1$, in the first case one has
\begin{equation}\label{eq:complete-double}
|\Sigma(r;G,\bfmu)|
\ll r^{1+\eps}(r,L_G(\bfmu))^{1/2},
\end{equation}
and in the second case one has
\begin{equation}\label{eq:complete-cube}
|\Sigma(r;G,\bfmu)|\ll
r^{3/2+\eps}\mathbf1_{r\mid L_G(\bfmu)}
\min\left\{r^{1/4},(r,M_G(\bfmu))\right\}.
\end{equation}
\end{enumerate}
The matrix $V$ and integer $A$ in \textup{(ii)} depend on $G$.
\end{theorem}

When $(r,6\Delta D_G(\bfmu))=1$, Theorem \ref{thm:complete}\textup{(i)} gives
square-root cancellation: $|\Sigma(r;G,\bfmu)|\ll r^{1+\eps}$. 
We prove this Theorem in
\cref{sec:complete}.

For comparison, consider the result of Loxton and Smith
\cite[Theorem~5]{Loxton1982}. For
$f(x,y)=G(x,y)+\mu_1x+\mu_2y$, let $\calD=D(\nabla f)$
denote the discriminant of the system $f_x=f_y=0$. 
When $\calD\ne0$, their bound becomes
\[
|\Sigma(p^k;G,\bfmu)|
\le4p^k(p^k,\calD^5)\qquad(k\ge2),
\]
for any prime $p$. This gives square-root cancellation
when $p\nmid\calD$, but does not
by itself give the dependence on the linear coefficients
needed for our averages. When $\bfmu=0$, one has $\calD=0$,
and the displayed right-hand side is $4p^{2k}$, which is
larger than the trivial bound $p^{2k}$.
In \cite[Section~7]{Loxton1982}, Loxton and Smith ask
for bounds covering the case $\calD=0$ and, when $\calD\ne0$,
for explicit dependence on the linear coefficients.
Under its coprimality conditions, \cref{thm:complete} answers
these questions in the binary cubic case, with explicit
dependence on both the linear coefficients and the cubic
discriminant.

The bounds in \cref{thm:complete} can be sharpened slightly in their
dependence on the linear coefficients and the cubic discriminant. We, 
however, find that the proofs of such refinements require more technical 
arguments and are not needed to prove \cref{thm:main}.
We plan to prove these stronger bounds in subsequent work.

\subsection{Notational conventions}\label{subsec:notation}

Throughout, $F$ and $\a$ are fixed and the integer $N\ge2$ tends to
infinity.  All statements involving $\eps$ are understood to hold
for every $\eps>0$.  Unless stated otherwise, implied constants
in Vinogradov's notation and $O(\cdot)$ may depend on
$\eps,\a$ and the coefficients of $F$.

We write $(a_1,\ldots,a_k)$ for the greatest common divisor of
the integers $a_1,\ldots,a_k$.  Vectors are denoted by bold
letters, and for $\bfx=(x_1,\ldots,x_n)$ we write
$|\bfx|=\max_i|x_i|$.  The notation $m\mid\bfx$ means that
$m$ divides every coordinate.  Sums over vectors in a box
$B\subseteq\R^n$ run over its integer points $B\cap\Z^n$.

For a prime $p$, let $\Z_p$ denote the ring of $p$-adic integers,
with its usual valuation $v_p$ and associated absolute value
$|t|_p=p^{-v_p(t)}$.  For vectors we write
$v_p(\bfx)=\min_i v_p(x_i)$.

\subsection*{Acknowledgements}

The author thanks D.~R.~Heath-Brown for a helpful conversation
about this work at the DECANT conference, held at ICMS in
Edinburgh in June 2025.

This work was supported in part by a grant of access to OpenAI
models through the ChatGPT for Academic Researchers program. The 
main ideas and the choice of techniques in this paper are
due to the author. ChatGPT-6.1 Sol was used to assist with
checking the manuscript for typographical and mathematical
errors and with revising the exposition.

\section{Preliminaries}\label{sec:preliminaries}

In the argument which follows, we will require rational approximations to $\a$ which
satisfy certain criteria. Namely, we would like the denominator $q = q_1q_2$ of our rational approximation
to $\a$ to contain a factor of a prescribed size, $q^{\eta} \le q_1 \le q^{\eta + \eps}$. This is not new and is in fact exactly
what Heath-Brown \cite{Heath2024} requires in his argument. Additionally, due to the form of Theorem
\ref{thm:complete}, we require that certain primes do not divide $q_2$. When $\a$ is a
quadratic irrational, we can guarantee the existence of such approximations.

\begin{lemma}\label{thm:smooth}
Let $\a$ be a real quadratic irrational and $\eta>0$.
For every sufficiently large $Q$ there are coprime integers $a,q$ with
\[
Q\ll_{\a,\eta}q\le Q,\qquad
|\a-a/q|\ll_{\a,\eta}(qQ)^{-1},
\]
such that $p^k\le q^\eta$ for every prime power
$p^k\parallel q$.
\end{lemma}

\begin{proof}
This is \cite[Theorem~3]{Heath2024}, with the prime-power
refinement noted immediately after \cite[Theorem~3]{Heath2010}.
\end{proof}

\begin{corollary}\label{cor:factorization}
Let $\mathcal P$ be a fixed finite set of primes, and let
$0<\beta<1$ and $\rho>0$ be fixed.  The approximations in
\cref{thm:smooth} can be chosen with a factorization $q=q_1q_2$
satisfying
\[
q^\beta\le q_1\le q^{\beta+\rho},
\qquad (q_2,\prod_{p\in\mathcal P}p)=1
\]
for all sufficiently large $Q$.
\end{corollary}

\begin{proof}
Choose $\eta<\min\{\rho,\beta/(|\mathcal P|+1)\}$ in
\cref{thm:smooth}.  Begin with
\[
q_1=\prod_{p\in\mathcal P}p^{v_p(q)}
\le q^{\eta|\mathcal P|}<q^\beta.
\]
Multiply $q_1$ by further prime factors of $q$, with multiplicity,
until it first reaches $q^\beta$.  Each factor is at most
$q^\eta$, so the resulting divisor satisfies
$q^\beta\le q_1\le q^{\beta+\eta}\le q^{\beta+\rho}$.
All prime-power factors supported on $\mathcal P$ lie in $q_1$,
which gives the required condition on $q_2=q/q_1$.
\end{proof}

\section{The two-dimensional \texorpdfstring{$A$}{A}-process}\label{sec:processes}

We begin with the two-dimensional $q$-analogue of the van der Corput $A$-process. 
This allows us to reduce the estimation of the binary quartic Weyl sum to that of 
a certain cubic Weyl sum. We first express $F(\bfm+\bft)-F(\bfm-\bft)$ in terms
of the partial derivatives of $F$.

For $\bft=(t_1,t_2)$ and $\bfm=(m_1,m_2)$, write
\[
G_{\bft}(\bfm)=t_1F_x(\bfm)+t_2F_y(\bfm) = \nabla F(\bfm) \cdot \bft.
\]
Thus $G_{\bft}$ is a binary cubic form in $\bfm$, whose coefficients
depend linearly on $\bft$.

For a quartic binary form $F$ note that 
\[F(\bfm \pm \bft) = F(\bfm) \pm \nabla F(\bfm) \cdot \bft + H(\bfm,\bft) \pm \nabla F(\bft) \cdot \bfm + F(\bft),\]
where $H(\bfm,\bft)$ is quadratic in both $\bfm$ and $\bft$. This follows from Taylor expansion about
both $\bfm$ and $\bft$. Hence, we have that
\begin{equation}\label{eq:quartic-difference}
F(\bfm+\bft)-F(\bfm-\bft)
=2\bigl(\nabla F(\bfm) \cdot \bft+\nabla F(\bft) \cdot \bfm\bigr) = 2(G_{\bft}(\bfm) + G_{\bfm}(\bft)).
\end{equation}

For $\bfh\in\Z^2$, a box $\calB'\subseteq\R^2$, and $a,q_1,q_2$ as below, set
\begin{equation}\label{eq:Th-def}
T_{\bfh}(\calB')=\sum_{\bfm\in\calB'}\e_{q_2}\Bigl(2a\bigl(G_{\bfh}(\bfm)+q_1^2\,\nabla F(\bfh)\cdot\bfm\bigr)\Bigr).
\end{equation}
\begin{lemma}\label{lem:A2}
Let $\a=a/q+\delta$ with $(a,q)=1$, $q=q_1q_2$ and $1\le q_1\le N/2$,
and put $H=\floor{N/(2q_1)}$.  Then
\[
|S_F(\a;N)|^2\ll(1+|\delta|N^4)^4q_1^2\Bigl(N^2+
\sum_{\substack{\bfh\in\Z^2\\0<|\bfh|<H}}\max_{\calB'}|T_{\bfh}(\calB')|\Bigr),
\]
where $\calB'$ runs over boxes contained in $(0,N]^2$.
\end{lemma}

\begin{proof}
For $0\le u,v\le N$ put
\[
S(u,v)=\sum_{\substack{0<m_1\le u\\0<m_2\le v}}\e_q(aF(\bfm)),
\qquad w(u,v)=\e(\delta F(u,v)).
\]
Riemann--Stieltjes integration by parts in each variable gives
\begin{align*}
S_F(\a;N)&=S(N,N)w(N,N)-\int_0^NS(u,N)\,\partial_1w(u,N)\dd u
-\int_0^NS(N,v)\,\partial_2w(N,v)\dd v\\
&\qquad+\int_0^N\!\!\int_0^NS(u,v)\,\partial_1\partial_2w(u,v)\dd u\dd v.
\end{align*}
On $[0,N]^2$ one has
\[
\partial_jw\ll|\delta|N^3,\qquad
\partial_1\partial_2w\ll|\delta|N^2+|\delta|^2N^6,
\]
so $|S_F(\a;N)|\ll(1+|\delta|N^4)^2\max_{u,v}|S(u,v)|$.

Fix $\calB'=(0,u]\times(0,v]$.  For $\bfh\in[1,H]^2$ and $\bfm\in\calB'$ put
$\bfn=\bfm-2q_1\bfh$; since $0<2q_1h_j\le2q_1H\le N$ we have
$\bfn\in(-N,N)^2$.  Hence
\[
H^2S(u,v)=\sum_{\bfn\in(-N,N)^2}\ \sum_{\substack{\bfh\in[1,H]^2\\ \bfn+2q_1\bfh\in\calB'}}
\e_q\bigl(aF(\bfn+2q_1\bfh)\bigr),
\]
and Cauchy's inequality gives
\[
H^4|S(u,v)|^2\le4N^2\sum_{\bfh,\bfh'\in[1,H]^2}\
\sum_{\substack{\bfn\\ \bfn+2q_1\bfh,\ \bfn+2q_1\bfh'\in\calB'}}
\e_q\Bigl(a\bigl(F(\bfn+2q_1\bfh)-F(\bfn+2q_1\bfh')\bigr)\Bigr).
\]
Put $\bfm=\bfn+q_1(\bfh+\bfh')$ and $\bfk=\bfh-\bfh'$, so that
$\bfn+2q_1\bfh=\bfm+q_1\bfk$ and $\bfn+2q_1\bfh'=\bfm-q_1\bfk$.  By
\eqref{eq:quartic-difference} and the linearity of $G_{\bft}$ in $\bft$,
\[
F(\bfm+q_1\bfk)-F(\bfm-q_1\bfk)=2\bigl(q_1G_{\bfk}(\bfm)+q_1^3G_{\bfm}(\bfk)\bigr)
=2q_1\bigl(G_{\bfk}(\bfm)+q_1^2\,\nabla F(\bfk)\cdot\bfm\bigr).
\]
Put $\calB'_{\bfk}=(\calB'-q_1\bfk)\cap(\calB'+q_1\bfk)$.
The summation conditions on $\bfn$ are then equivalent to $\bfm\in\calB'_{\bfk}$, so by
\eqref{eq:Th-def} the inner sum over $\bfn$ is
\[
\sum_{\bfm\in\calB'_{\bfk}}
\e_q\bigl(a(F(\bfm+q_1\bfk)-F(\bfm-q_1\bfk))\bigr)
=T_{\bfk}(\calB'_{\bfk}).
\]
The box $\calB'_{\bfk}$ lies in $\calB'$, since $\bfm$ is the
midpoint of $\bfm\pm q_1\bfk$.  Each $\bfk$ with $|\bfk|<H$
arises from at most $H^2$ pairs $(\bfh,\bfh')$, and
$T_{\bfzero}(\calB'_{\bfzero})=\#(\calB'\cap\Z^2)\le N^2$.
Therefore
\[
\begin{aligned}
H^4|S(u,v)|^2
&\le4N^2H^2\sum_{|\bfk|<H}|T_{\bfk}(\calB'_{\bfk})|\\
&\le4N^2H^2\Bigl(N^2+
\sum_{0<|\bfk|<H}\max_{\calB''}|T_{\bfk}(\calB'')|\Bigr).
\end{aligned}
\]
Since $q_1\asymp N/H$, dividing by $H^4$ and renaming $\bfk$
as $\bfh$ gives the stated bound.
\end{proof}

It remains to estimate the cubic sums $T_{\bfh}(\calB')$.

\section{The two-dimensional \texorpdfstring{$B$}{B}-process}\label{sec:completion}\label{subsec:modulus}

We now estimate the cubic Weyl sums $T_{\bfh}(\calB')$.
We begin by bounding them in terms of complete Weyl sums.

For $r\ge1$ put
\begin{equation}\label{eq:weights}
w_r(t)=\min\Bigl(N,\frac{2r}{|t|}\Bigr)\quad(t\ne0),\qquad w_r(0)=N,\qquad
W_r(\bfj)=w_r(j_1)w_r(j_2).
\end{equation}
Write $\mathcal J_r=((-r/2,r/2]\cap\Z)^2$.

\begin{lemma}\label{lem:B2}
Let $G$ be an integral binary cubic form, let $\bflam\in\Z^2$,
and let $\calB'\subseteq(0,N]^2$ be a box.  Then
\begin{equation}\label{eq:B2}
\left|\sum_{\bfm\in\calB'}\e_r\bigl(G(\bfm)+\bflam\cdot\bfm\bigr)\right|
\le\frac1{r^2}\sum_{\bfj\in\mathcal J_r}
W_r(\bfj)\,|\Sigma(r;G,\bflam+\bfj)|.
\end{equation}
\end{lemma}

\begin{proof}
Grouping the summands according to $\bfm\bmod r$ and using
orthogonality, we obtain
\[
\begin{aligned}
\sum_{\bfm\in\calB'}\e_r\bigl(G(\bfm)+\bflam\cdot\bfm\bigr) &=\sum_{\bfm'\bmod r}
\e_r\bigl(G(\bfm')+\bflam\cdot\bfm'\bigr)
\sum_{\substack{\bfm\in\calB'\\\bfm\equiv\bfm'\pmod r}}1\\
&=\frac1{r^2}\sum_{\bfj\in\mathcal J_r}
\Sigma(r;G,\bflam+\bfj)
\sum_{\bfm\in\calB'}\e_r(-\bfj\cdot\bfm).
\end{aligned}
\]
The inner geometric sums have absolute value at most $W_r(\bfj)$.
\end{proof}

For an integral cubic $G$ and $\bflam\in\Z^2$, write
\begin{equation}\label{eq:Psi-def}
\Psi_N(r;G,\bflam)=\frac1{r^2}\sum_{\bfj\in\mathcal J_r}
W_r(\bfj)\,|\Sigma(r;G,\bflam+\bfj)|.
\end{equation}
In particular, \cref{lem:B2} applied to \eqref{eq:Th-def} gives
\[
\max_{\calB'}|T_{\bfh}(\calB')|
\le\Psi_N\bigl(q_2;2aG_{\bfh},2aq_1^2\nabla F(\bfh)\bigr).
\]
Before summing over $\bfh$, we estimate $\Psi_N$ for a fixed
cubic. To apply \cref{thm:complete}\textup{(i)}, we first
count the pairs $\bfmu$ in a translated rectangle
satisfying $m\mid\bfmu$ and $g\mid D_G(\bfmu)$.

\begin{lemma}\label{lem:quadratic-count}
Let $D$ be an integral binary quadratic form whose coefficients
have greatest common divisor one, with
$\disc(D)\ne0$.  Let $g$ be a positive odd integer, let $m$ be
a positive integer, and let $\mathcal R=I_1\times I_2$, where $I_i$
is a real interval of length $M_i\ge0$.  Then
\[
\#\{\bfmu\in\mathcal R\cap\Z^2:
m\mid\bfmu,\ g\mid D(\bfmu)\}
\ll g^\eps
\left(\frac{M_1M_2}{m^2\bigl(g/(g,m^2)\bigr)^{1/2}}
+\frac{M_1+M_2}{m}+1\right).
\]
If $(g,\disc(D))=1$, the exponent $1/2$ on $g/(g,m^2)$
may be replaced by $1$.
\end{lemma}

\begin{proof}
The case $g=1$ is trivial.  Since $D$ is quadratic, we have
\[
g\mid D(m\boldsymbol\nu)
\quad\Longleftrightarrow\quad
\frac{g}{(g,m^2)}\mid D(\boldsymbol\nu).
\]
Thus writing $\bfmu=m\boldsymbol\nu$ replaces $g$ by $g/(g,m^2)$
and the interval lengths $M_i$ by $M_i/m$.  It suffices to prove
the bound for $m=1$.

Suppose henceforth that $m=1$ and $g>1$.  We adapt the proof of
\cite[Lemma~2.4]{Bruedern2014}.
Write $D(X,Y)=AX^2+BXY+CY^2$.
To complete the square modulo $g$, we first choose new variables
$S,T$ so that the coefficient of $S^2$ is relatively prime to $g$.
Since the coefficients of $D$ have greatest common divisor one,
$D(1,T)$ is a nonzero polynomial modulo each prime $p\mid g$.
It has at most two roots, and $p\ge3$, so there is a residue
at which it does not vanish.  The Chinese remainder theorem
therefore gives an integer $c$ with $(D(1,c),g)=1$.
Make the substitution $S=X$, $T=Y-cX$.  Then
\[
D(X,Y)=A_1S^2+B_1ST+CT^2,
\]
where $A_1=D(1,c)$, $B_1=B+2Cc$, and $(A_1,g)=1$.

Since $B_1^2-4A_1C=B^2-4AC=\disc(D)$, completing the square gives
\[
4A_1 D(X,Y)=(2A_1S+B_1T)^2-\disc(D)T^2.
\]
As $(2A_1,g)=1$, choose an integer $\tau$ with
$2A_1\tau\equiv B_1\pmod g$ and put
\[
U=S+\tau T,\qquad V=T.
\]
The congruence $D(X,Y)\equiv0\pmod g$ is then equivalent to
\[
U^2\equiv bV^2\pmod g, \quad \text{where} \quad
4A_1^2 b\equiv\disc(D)\pmod g.\]

For a solution $U,V$ to the above congruence, put $h=(U,V,g)$ and write
\[
U=hu,\qquad V=hv,\qquad q=\frac{g}{(g,h^2)}.
\]
Then $u^2\equiv bv^2\pmod q$ and $(u,v,q)=1$, whence
$(v,q)=1$.  Thus we may write $u\equiv zv\pmod q$, where
\begin{equation}\label{eq:quadratic-root-congruence}
z^2\equiv b\pmod q.
\end{equation}
We first show that the roots of \eqref{eq:quadratic-root-congruence}
lie in $O(q^\eps)$
residue classes modulo a divisor $d$ of $q$, where $d\ge q^{1/2}$.
We then count the corresponding pairs in the original rectangle
$\mathcal R$ for each of these residue classes.

If $b$ is not a square modulo $q$, then there are no solutions to \eqref{eq:quadratic-root-congruence}.
  Otherwise choose one root $z_0$.  Every root $z$ satisfies
\[
q\mid(z-z_0)(z+z_0).
\]
For each $p^\ell\parallel q$, at least one factor is divisible
by $p^{\lceil\ell/2\rceil}$.  Hence every root $z$ satisfies
\[
z\equiv z_0\quad\text{or}\quad z\equiv-z_0
\pmod{p^{\lceil\ell/2\rceil}}.
\]
Put
\[
d=\prod_{p^\ell\parallel q}p^{\lceil\ell/2\rceil}
\ge q^{1/2}.
\]
Each choice of signs in these congruences determines, by the Chinese
remainder theorem, one residue class modulo $d$.  Thus, by the divisor
bound, the roots of
\eqref{eq:quadratic-root-congruence} lie in $O(q^\eps)$
residue classes modulo $d$.

For each such residue $z$, the corresponding pairs $U,V$ satisfy
\[
h\mid V,\qquad U\equiv zV\pmod{hd}.
\]
Modulo $hd$, there are $d$ choices for $V$, and each determines
$U$ uniquely.  Returning to the original variables $X,Y$, this
system becomes
\[
h\mid Y-cX,\qquad X+(\tau-z)(Y-cX)\equiv0\pmod{hd}.
\]
Both substitutions are bijective modulo $hd$, so this system
also has exactly $d$ solutions modulo $hd$.

We count these solutions in the original rectangle $\mathcal R$,
choosing $Y$ first and then $X$ for each fixed $Y$.
The following argument applies to any system of homogeneous linear
congruences in $X,Y$.

Sums and differences of solutions are again solutions, so the
possible values of $Y$ are closed under addition and subtraction.
Let $d_2$ be the smallest positive possible value of $Y$.
It exists since $(0,hd)$ is a solution.
Thus the possible values of $Y$ are exactly the multiples of $d_2$.

For each $Y\in d_2\Z$, choose a solution $(X_0(Y),Y)$.
Another pair $(X,Y)$ is a solution precisely when
$(X-X_0(Y),0)$ is a solution.
Let $d_1$ be the smallest positive integer for which $(d_1,0)$
is a solution.  For this fixed $Y$, the pair $(X,Y)$ is a
solution precisely when
\[
X\equiv X_0(Y)\pmod{d_1}.
\]
The modulus $d_1$ is independent of $Y$.

Since $(hd,0)$ and $(0,hd)$ are solutions, both $d_1$ and $d_2$
divide $hd$.  Modulo $hd$, there are $hd/d_2$ possible values of $Y$
and $hd/d_1$ values of $X$ for each.  We already counted exactly
$d$ solutions modulo $hd$, so
\[
d=\frac{hd}{d_2}\frac{hd}{d_1},\qquad d_1d_2=h^2d.
\]
In $I_2$ there are at most $1+M_2/d_2$ possible values of $Y$.
For each, there are at most $1+M_1/d_1$ values of $X$ in $I_1$.
Consequently,
\begin{equation}\label{eq:linear-congruence-count}
\begin{aligned}
\#\{(X,Y)\in\mathcal R\cap\Z^2:
d_2\mid Y,\ X\equiv X_0(Y)\pmod{d_1}\}
&\le\left(\frac{M_1}{d_1}+1\right)
     \left(\frac{M_2}{d_2}+1\right)
\\
&\le\frac{M_1M_2}{d_1d_2}+M_1+M_2+1.
\end{aligned}
\end{equation}
Since $d\ge q^{1/2}$ and $(g,h^2)\le h^2$, we have
\[
h^2d\ge h^2\left(\frac{g}{(g,h^2)}\right)^{1/2}
\ge hg^{1/2}\ge g^{1/2}.
\]
By the divisor bound, there are $O(g^\eps)$ such systems of
congruences. Using $d_1d_2=h^2d$ in
\eqref{eq:linear-congruence-count} and summing over them proves
the first bound.
Finally, if $(g,\disc(D))=1$, then $(q,b)=1$.  Hence $p$
cannot divide both $z-z_0$ and $z+z_0$, so one factor is
divisible by $p^\ell$.  We may therefore take $d=q$, so
$h^2d\ge g$ and the first term in
\eqref{eq:linear-congruence-count} is at most $M_1M_2/g$.
\end{proof}

We now combine this count with \cref{thm:complete} to bound $\Psi_N$.

\begin{lemma}\label{lem:completion}
Let $G$ be a nonzero integral binary cubic form, let
$\bflam\in\Z^2$, and let $r\ge1$ be an integer and $N\ge2$.
Write $\Delta=\disc(G)$.
\begin{enumerate}[label=\textup{(\roman*)},leftmargin=*]
\item If $\Delta\ne0$ and $(r,6)=1$, then
\begin{equation}\label{eq:Psi-bound}
\Psi_N(r;G,\bflam)
\ll r^\eps(r^4,\Delta)^{1/4}
\bigl(r+Nr^{1/4}+N^2r^{-1/2}\bigr).
\end{equation}
\item If $\Delta=0$, choose $V$ and $A$ as in
\cref{thm:complete}\textup{(ii)}.  For $(r,6A\det V)=1$, one has
\[
\Psi_N(r;G,\bflam)\ll r^\eps
\begin{cases}
r+Nr^{1/2}+N^2r^{-1/2}, &G\circ V=AX^2Y,\\
Nr^{1/2}+N^2r^{-1/4}, &G\circ V=AX^3.
\end{cases}
\]
\end{enumerate}
\end{lemma}

\begin{proof}
The case $r=1$ is immediate.  Suppose first that $\Delta\ne0$,
and put $D=D_G$ and $B=(r^4,\Delta)^{1/4}$.
The coefficients of $D$ have greatest common divisor one, and
$\disc(D)=-3\Delta/c_G^2\ne0$ by \eqref{eq:covariants}.
For $m\mid r$, $g\mid r$, put
$\delta(m,g)=\min(r^{1/4}m^{1/2},g^{1/2})$.
By \cref{thm:complete}\textup{(i)} and \eqref{eq:Psi-def},
\[
\begin{aligned}
\Psi_N(r;G,\bflam)
&\ll Br^{\eps-1}
\sum_{\bfj\in\mathcal J_r}W_r(\bfj)
\delta\bigl((r,\bflam+\bfj),(r,D(\bflam+\bfj))\bigr)\\
&\le Br^{\eps-1}
\sum_{\bfj\in\mathcal J_r}W_r(\bfj)
\sum_{\substack{m\mid r,\ m\mid\bflam+\bfj\\
g\mid r,\ g\mid D(\bflam+\bfj)}}\delta(m,g)\\
&=Br^{\eps-1}
\sum_{\substack{m\mid r\\g\mid r}}\delta(m,g)
\sum_{\substack{\bfj\in\mathcal J_r\\
m\mid\bflam+\bfj,\ g\mid D(\bflam+\bfj)}}W_r(\bfj).
\end{aligned}
\]
Split the inner sum into $O(\log(2r)^2)$ dyadic rectangles.
On each rectangle let $M_i$ be the interval lengths and $w_i$
upper bounds for $w_r(j_i)$.  Then $w_i\le N$ and $w_iM_i\ll r$,
including the zero interval, where $M_i=0$.
By \cref{lem:quadratic-count}, the contribution of this rectangle,
including $\delta(m,g)$, is at most $O(g^\eps)$ times
\[
\delta(m,g)\left(
\frac{r^2(g,m^2)^{1/2}}{m^2g^{1/2}}+\frac{Nr}{m}+N^2\right)
\ll r^2+Nr^{5/4}+N^2r^{1/2}.
\]
Summing over the divisors and dyadic rectangles and using the
divisor bound gives \eqref{eq:Psi-bound}.

Suppose now that $\Delta=0$, and write
$V^{\mathsf T}\bfmu=(M_G(\bfmu),L_G(\bfmu))^{\mathsf T}$.
The map
$\bfmu\mapsto(M_G(\bfmu),L_G(\bfmu))$ is then invertible
modulo every divisor of $r$.
We apply \eqref{eq:complete-double} or \eqref{eq:complete-cube}
and use the same dyadic rectangles as above.

In the double-root case we use
\[
(r,L_G(\bfmu))^{1/2}
\le\sum_{\substack{g\mid r\\g\mid L_G(\bfmu)}}g^{1/2}.
\]
Since the substitution is invertible modulo $g$, there are
exactly $g$ pairs $\bfmu\bmod g$ satisfying
$L_G(\bfmu)\equiv0\pmod g$.
The counting argument in the proof of \cref{lem:quadratic-count}
therefore bounds the weighted number of pairs on each dyadic
rectangle by $O(r^2/g+Nr+N^2)$.
Thus, by the divisor bound,
\[
\begin{aligned}
\Psi_N(r;G,\bflam)
&\ll r^{\eps-1}
\sum_{g\mid r}g^{1/2}\left(\frac{r^2}g+Nr+N^2\right)\\
&\ll r^\eps\bigl(r+Nr^{1/2}+N^2r^{-1/2}\bigr).
\end{aligned}
\]

In the cube case, \eqref{eq:complete-cube} restricts the sum
to $r\mid L_G(\bfmu)$.  For each $g\mid r$, imposing also
$g\mid M_G(\bfmu)$ gives $g\mid\bfmu$ by the inverse
substitution, since $g\mid L_G(\bfmu)$ as well.
Writing $\bfmu=g\boldsymbol\nu$
leaves the single congruence $(r/g)\mid L_G(\boldsymbol\nu)$,
which has $r/g$ solutions modulo $r/g$.
On a dyadic rectangle $\mathcal R$ with interval lengths $M_1,M_2$, the
relation $\bfmu=\bflam+\bfj=g\boldsymbol\nu$ places
$\nu_i$ in an interval of length $M_i/g$.
The argument for counting solutions of linear congruences in the
proof of \cref{lem:quadratic-count} therefore gives
$O(M_1M_2/(rg)+(M_1+M_2)/g+1)$ such pairs.
Multiplying by $w_1w_2$ and using $w_iM_i\ll r$ and $w_i\le N$
gives
\[
\begin{aligned}
\sum_{\substack{\bfj\in\mathcal R\cap\mathcal J_r\\
g\mid\bflam+\bfj,\ r\mid L_G(\bflam+\bfj)}}W_r(\bfj)
&\ll w_1w_2\left(\frac{M_1M_2}{rg}
+\frac{M_1+M_2}{g}+1\right) \ll \frac{Nr}{g}+N^2.
\end{aligned}
\]
Putting $\delta_g=\min(r^{1/4},g)$ and using the divisor bound, we obtain
\[
\begin{aligned}
\Psi_N(r;G,\bflam)
&\ll r^{\eps-1/2}
\sum_{g\mid r}\delta_g
\left(\frac{Nr}g+N^2\right)\\
&\ll r^\eps
\bigl(Nr^{1/2}+N^2r^{-1/4}\bigr).
\end{aligned}
\]
\end{proof}

Since $F$ is non-degenerate, each $G_{\bfh}$ with $\bfh\ne0$
is a nonzero cubic. Otherwise a nonsingular change of variables
would make $F$ independent of one variable, forcing it to be a
fourth power.

We first treat quartics for which every linear factor over
$\mathbb C$ occurs at most twice.  To average the bounds in
\cref{lem:completion} over $\bfh$, we write
\[
\sum_{0<|\bfh|<H}\max_{\calB'}|T_{\bfh}(\calB')|=S_1+S_2,
\]
where
\[
\begin{aligned}
S_1&=\sum_{\substack{0<|\bfh|<H\\\disc_{\bfm}(G_{\bfh})=0}}
\max_{\calB'}|T_{\bfh}(\calB')|,\\
S_2&=\sum_{\substack{0<|\bfh|<H\\\disc_{\bfm}(G_{\bfh})\ne0}}
\max_{\calB'}|T_{\bfh}(\calB')|.
\end{aligned}
\]
Thus \cref{lem:A2} gives
\begin{equation}\label{eq:reduction}
|S_F(\a;N)|^2
\ll(1+|\delta|N^4)^4q_1^2(N^2+S_1+S_2).
\end{equation}

To estimate $S_1$ and $S_2$, we must understand the discriminant
of $G_{\bfh}$ as a function of $\bfh$.  Put
\[
D_F(\bfh)=\disc_{\bfm}(G_{\bfh}).
\]
For $S_1$ we need to count the shifts for which $D_F(\bfh)=0$;
for $S_2$ we must average the factor $(q_2^4,D_F(\bfh))^{1/4}$.
The following Lemma shows that $D_F$ is nonzero and that each
of its linear factors occurs at most twice.
The first property confines the shifts in $S_1$ to finitely many
rational lines; the second controls the congruence conditions
needed to estimate $S_2$.

\begin{lemma}\label{lem:polar-discriminant}
Suppose that every linear factor of $F$ over $\mathbb C$
occurs at most twice.  The nonzero integer shifts $\bfh$
for which $G_{\bfh}$ has a repeated linear factor lie on at most
four rational lines through the origin.

Moreover,
\[
D_F(\bfh)=\disc_{\bfm}(G_{\bfh})
\]
is a nonzero homogeneous polynomial of degree four, and each
linear factor in its factorization over $\mathbb C$ occurs at
most twice.
\end{lemma}

\begin{proof}
The cubic discriminant is homogeneous of degree four in the
coefficients.  Since the coefficients of $G_{\bfh}$ are linear
in $\bfh$, it follows that $D_F$ is either zero or homogeneous
of degree four in $\bfh$.

If $\widetilde F=F\circ M$ for a nonsingular matrix $M$, the
chain rule and the transformation formula for the cubic
discriminant give
\[
\bft\cdot\nabla\widetilde F=G_{M\bft}\circ M,
\qquad
D_{\widetilde F}(\bft)=(\det M)^6D_F(M\bft).
\]
Thus a nonsingular change of variables preserves the
multiplicities of the linear factors of $D_F$.

Suppose first that $\disc(F)\ne0$.
The forms $F_x,F_y$ have no common nonzero zero, since Euler's
identity $xF_x+yF_y=4F$ would otherwise give a repeated factor of $F$.
Writing $t=x/y$ and $z=h_1/h_2$ for $h_2\ne0$, set
\[
f_z(t)=\frac{G_{\bfh}(t,1)}{h_2}
      =F_y(t,1)+zF_x(t,1).
\]
For a fixed $z_0$, we choose a nonsingular linear change in
$x,y$, with the corresponding change in $\bfh$, so that
$h_2\ne0$ and the coefficient of $x^3$ in $G_{\bfh}(x,y)$
is nonzero.
Thus $f_{z_0}(t)$ has degree three.
If $f_{z_0}$ has no repeated root, then $D_F$ is nonzero at
this ratio.  Otherwise let $t_0$ be a repeated root, so that
$F_x(t_0,1)\ne0$.  Write
\[
f_{z_0+u}(t_0+t)=t^2(at+b)+uF_x(t_0+t,1),\qquad a\ne0.
\]
Here $b\ne0$ for a double root and $b=0$ for a triple root.
Substituting these coefficients into the cubic discriminant gives
\[
\disc(f_{z_0+u})=
\begin{cases}
-4b^3F_x(t_0,1)u+O(u^2),&t_0\text{ is a double root},\\
-27a^2F_x(t_0,1)^2u^2+O(u^3),&t_0\text{ is a triple root}.
\end{cases}
\]
The nonzero leading coefficients show that $D_F$ is not
identically zero and that the order of vanishing is respectively
one or two.  The transformation formula above gives the same
orders for $D_F$.

Suppose now that $\disc(F)=0$.  If $F$ has one double factor
and two simple factors, a nonsingular change of variables
over $\mathbb C$ puts it in the form
\[
F=x^2(Ax^2+Bxy+Cy^2),\qquad C(B^2-4AC)\ne0.
\]
Writing $\bfh=(u,v)$, direct calculation gives
\[
D_F(u,v)=4C^2u^2
\bigl((9B^2-32AC)u^2+4BCuv+4C^2v^2\bigr).
\]
The quadratic factor is not divisible by $u$, and its
discriminant is $128C^2(4AC-B^2)\ne0$.
If $F$ has two double factors, write $F=cQ^2$, where $Q$ is
a quadratic form with $\disc(Q)\ne0$.  Using
$G_{\bfh}=2cQ(\bfh\cdot\nabla Q)$ gives
\[
D_F(\bfh)=16c^4\disc(Q)^3Q(\bfh)^2.
\]
In both cases $D_F$ is nonzero and every linear factor occurs
at most twice.  Finally, its degree is four, so its nonzero
integer zeros lie on at most four rational lines through the
origin.
\end{proof}

We now fix the primes to be excluded from $q_2$.  We exclude
$2$ and $3$ to satisfy the coprimality conditions in
\cref{thm:complete}.

For $S_1$, choose a primitive integer vector $\bfv$ on each
rational line on which $D_F$ vanishes.  By
\cref{thm:complete}\textup{(ii)}, we can choose an integer
matrix $V_{\bfv}$ with nonzero determinant and a nonzero
integer $A_{\bfv}$ such that
\[
G_{\bfv}\circ V_{\bfv}=A_{\bfv}X^2Y
\quad\text{or}\quad A_{\bfv}X^3.
\]
We exclude the primes dividing $A_{\bfv}\det V_{\bfv}$ for
each of these finitely many choices.  The substitutions are
then invertible modulo $q_2$, and each $A_{\bfv}$ is relatively
prime to $q_2$, as required in that Theorem.

For $S_2$, we will count shifts satisfying $g\mid D_F(\bfh)$.
Factor $D_F$ over $\Q$ as
\[
D_F=c\prod_{\ell=1}^s P_\ell^{m_\ell},
\qquad c\in\Z\setminus\{0\},\quad m_\ell\in\{1,2\}.
\]
Here the $P_\ell$ are distinct irreducible homogeneous
polynomials with integer coefficients of gcd one.
Put $R_F=\prod_\ell P_\ell$.  Each factor occurs once in
$R_F$, so $\disc(R_F)\ne0$.  We exclude primes dividing $c$,
since these divide $D_F(\bfh)$ for every shift.  We also exclude
primes dividing $\disc(R_F)$, so that the roots of $R_F$ modulo
each remaining prime are simple.  Hensel's Lemma then applies
in the congruence count below.

Thus let $\mathcal P_F$ be the set of prime divisors of
\[
6c\disc(R_F)\prod_{\bfv}A_{\bfv}\det V_{\bfv},
\]
where the product takes one chosen vector on each rational line
above.  This is a fixed finite set depending only on $F$.
By \cref{cor:factorization}, we may place every prime-power
factor of $q$ belonging to $\mathcal P_F$ in $q_1$, so that
$q_2$ has no prime factor in $\mathcal P_F$.
The remaining divisibility conditions $p^j\mid D_F(\bfh)$,
for $p\mid q_2$, will be handled by the average over shifts
in $S_2$.

We estimate $S_1$ and $S_2$ separately.

\begin{lemma}\label{lem:singular-average}
Suppose every linear factor of $F$ over $\mathbb C$ occurs at
most twice.  If $q_2$ has no prime factor in $\mathcal P_F$, then
\[
S_1\ll Hq_2^\eps
\bigl(q_2+Nq_2^{1/2}+N^2q_2^{-1/4}\bigr).
\]
\end{lemma}

\begin{proof}
By \cref{lem:polar-discriminant}, we have
\[
S_1=\sum_{\bfv}\sum_{0<|t|<H/|\bfv|}
\max_{\calB'}|T_{t\bfv}(\calB')|,
\]
where the outer sum runs over the at most four primitive
vectors chosen above.  Fix $\bfv$ and put $d=(q_2,t)$ and $r=q_2/d$.
After cancelling $d$ in \eqref{eq:Th-def}, the cubic part is
$(2at/d)G_{\bfv}$ modulo $r$, with $(2at/d,r)=1$.
Thus \cref{lem:B2} and \cref{lem:completion}\textup{(ii)} give
\[
\begin{aligned}
\max_{\calB'}|T_{t\bfv}(\calB')|
&\ll r^\eps\bigl(r+Nr^{1/2}+N^2r^{-1/4}\bigr)\\
&\ll q_2^\eps\bigl(q_2+Nq_2^{1/2}
     +N^2q_2^{-1/4}d^{1/4}\bigr).
\end{aligned}
\]
The divisor bound gives
\[
\sum_{1\le t<H}(q_2,t)^{1/4}
\le H\sum_{d\mid q_2}d^{-3/4}
\ll Hq_2^\eps.
\]
Summing over $t$ and the finitely many choices of $\bfv$ proves
the result.
\end{proof}

\begin{lemma}\label{lem:nonsingular-average}
Suppose every linear factor of $F$ over $\mathbb C$ occurs at
most twice.  If $q_2$ has no prime factor in $\mathcal P_F$, then
\begin{equation}\label{eq:shift-average}
S_2\ll H^2q_2^\eps
\bigl(q_2+Nq_2^{1/4}+N^2q_2^{-1/2}\bigr).
\end{equation}
\end{lemma}

\begin{proof}
For $D_F(\bfh)\ne0$, apply \cref{lem:B2} and
\cref{lem:completion}\textup{(i)} at modulus $q_2$ to $2aG_{\bfh}$.
Since $(2a,q_2)=1$ and $\disc(2aG_{\bfh})=(2a)^4D_F(\bfh)$,
it suffices to prove
\[
\sum_{\substack{0<|\bfh|<H\\D_F(\bfh)\ne0}}
(q_2^4,D_F(\bfh))^{1/4}\ll H^2q_2^\eps.
\]

We begin by counting the pairs satisfying $g\mid D_F(\bfh)$.
Let $g\ge1$ have no prime factor in $\mathcal P_F$, and let
$\mathcal R=I_1\times I_2$ have interval lengths $M_1,M_2\ge0$.
We will use the bound
\[
\#\{\bfh\in\mathcal R\cap\Z^2:g\mid D_F(\bfh)\}
\ll g^\eps
\left(\frac{M_1M_2}{g^{1/2}}+M_1+M_2+1\right).
\]
We use the counting argument in the proof of \cref{lem:quadratic-count}.
For each $p^j\parallel g$, it suffices to show that every solution of
$p^j\mid D_F(\bfh)$ satisfies one of $O(j+1)$ systems of homogeneous
linear congruences, each having at most $p^{3j/2}$ solutions modulo $p^j$.

The condition $p^{\lceil j/4\rceil}\mid\bfh$ accounts for shifts
of large common valuation.  It has
$p^{2(j-\lceil j/4\rceil)}\le p^{3j/2}$ solutions modulo $p^j$.
Otherwise write $\bfh=p^s\mathbf v$ with $p\nmid(v_1,v_2)$
and $4s<j$.  Suppose first that $p\nmid v_2$, so that we may
use $z=v_1/v_2$ modulo powers of $p$.
The congruence $D_F(z,1)\equiv0\pmod p$ forces $z$
to be a root of $R_F(z,1)$ modulo $p$.
Since $p\nmid\disc(R_F)$, this root is simple and lifts to
$z_0\in\Z_p$ by Hensel's Lemma.  The corresponding factor
occurs at most twice in $D_F$, and all other factors are
relatively prime to $p$ at $\mathbf v$.  Hence
\[
v_p(D_F(\mathbf v))\le2v_p(v_1-z_0v_2).
\]
If $p\mid v_2$, then $p\nmid v_1$, and interchanging the
coordinates gives the same conclusion with
$L(\mathbf v)=v_2-z_0v_1$.
In the first case put $L(\mathbf v)=v_1-z_0v_2$.
In either case a coefficient of $L$ is relatively prime to $p$.
Put $t=\lceil(j-4s)/2\rceil$.  The corresponding shifts satisfy
\[
p^s\mid\bfh,\qquad
L(\bfh)\equiv0\pmod{p^{s+t}}.
\]
Since a coefficient of $L$ is relatively prime to $p$, these congruences have
\[
p^{j-s}p^{j-s-t}=p^{2j-2s-t}\le p^{3j/2}
\]
solutions modulo $p^j$.
There are at most four roots to consider and $O(j+1)$ choices of $s$.
The Chinese remainder theorem, the divisor bound and
\eqref{eq:linear-congruence-count}, as in the proof of
\cref{lem:quadratic-count}, give the preceding bound.

For $X>1$ and a positive integer $R$ having no prime factor in
$\mathcal P_F$, we now average $(R,D_F(\bfh))^{1/4}$ over
$0<|\bfh|<X$ with $D_F(\bfh)\ne0$.
Every divisor $g$ of both $R$ and $D_F(\bfh)$ satisfies
$g\ll_F X^4$.  The preceding count therefore gives
\begin{equation}\label{eq:discriminant-average}
\begin{aligned}
\sum_{\substack{0<|\bfh|<X\\D_F(\bfh)\ne0}}
(R,D_F(\bfh))^{1/4}
&\ll\sum_{\substack{g\mid R\\g\ll_F X^4}}
g^{1/4+\eps}\left(X^2g^{-1/2}+X+1\right)\\
&\ll X^2R^\eps,
\end{aligned}
\end{equation}
where we used $g^{1/4}\ll_F X$ and the divisor bound.

Taking $R=q_2^4$ and $X=H$ in
\eqref{eq:discriminant-average} proves the required estimate.
If $H=1$, the shift sum is empty.
\end{proof}

It remains to treat quartics with a triple linear factor.
That factor is defined over $\Q$, since distinct conjugate
factors of multiplicity three would have total degree at least
six.  Thus, unless $F$ is a fourth power, we can write
$F=cL^3M$, where $L,M$ are independent primitive integral
linear forms and $c$ is a nonzero integer.
Every $G_{\bfh}$ then has a repeated factor, so we use
\cref{lem:mixed-cubic} to estimate the sum over shifts.

\begin{lemma}\label{lem:triple-factor-average}
Suppose that $F=cL^3M$ as above.  There is a fixed finite set
$\mathcal P_F$ of primes such that, if $q_2$ has no prime
factor in $\mathcal P_F$, then
\[
\begin{aligned}
\sum_{0<|\bfh|<H}\max_{\calB'}|T_{\bfh}(\calB')|
\ll q_2^\eps\bigl(&H^2(q_2+N+N^2q_2^{-1/2})\\
&+H(q_2+Nq_2^{1/2}+N^2q_2^{-1/4})\bigr).
\end{aligned}
\]
\end{lemma}

\begin{proof}
The result is immediate if $H=1$ or $q_2=1$.
To apply \cref{lem:mixed-cubic}, write
$L(x,y)=\ell_1x+\ell_2y$ and $M(x,y)=m_1x+m_2y$, and put
\[
D=\ell_1m_2-\ell_2m_1\ne0,\qquad
V=\begin{pmatrix}m_2&-\ell_2\\-m_1&\ell_1\end{pmatrix}.
\]
Then $L(V(X,Y)^{\mathsf T})=DX$,
$M(V(X,Y)^{\mathsf T})=DY$, and $\det V=D$.  Hence
\[
F\circ V=cD^4X^3Y,\qquad
G_{\bfh}\circ V=cD^3\bigl(M(\bfh)X^3+3L(\bfh)X^2Y\bigr).
\]
Let $\mathcal P_F$ consist of the prime divisors of $6cD$.
Put $r=q_2$ and suppose first that $L(\bfh)\ne0$.
The substitution $\bfs=V(X,Y)^{\mathsf T}$ is bijective
modulo $r$.  For $\bfmu\in\Z^2$, write
$\boldsymbol\nu=V^{\mathsf T}\bfmu$.  Since $(2a,r)=1$,
the proof of \cref{lem:mixed-cubic} gives
\[
|\Sigma(r;2aG_{\bfh},\bfmu)|
\ll r^{1+\eps}(r,L(\bfh))
\sum_{\substack{m^2\mid r\\m\mid\nu_1,\ m^2\mid\nu_2}}m.
\]

To estimate $\Psi_N$, we count the original coefficient
pairs $\bfmu$ on the translated rectangles.
For $m^2\mid r$, the conditions $m\mid\nu_1$ and
$m^2\mid\nu_2$ imply $\bfmu=m\boldsymbol\eta$, since
$(\det V,m)=1$.  They are then equivalent to the single
congruence
\[
m\mid(V^{\mathsf T}\boldsymbol\eta)_2.
\]
The coefficients are not both divisible by any prime dividing
$m$, so there are exactly $m$ solutions modulo $m$.
On a translated rectangle $\mathcal R$ in $\bfmu$ with
interval lengths $M_1,M_2$, the count
\eqref{eq:linear-congruence-count}, applied in
$\boldsymbol\eta$, therefore gives
\[
\#\{\bfmu\in\mathcal R\cap\Z^2:
m\mid\nu_1,\ m^2\mid\nu_2\}
\ll \frac{M_1M_2}{m^3}+\frac{M_1+M_2}{m}+1.
\]
On each translated dyadic rectangle $\mathcal R$ arising from
$\bfmu=\bflam+\bfj$, we have $w_r(\mu_i-\lambda_i)\le N$ and
$M_iw_r(\mu_i-\lambda_i)\ll r$.  These bounds give
\[
\sum_{\substack{\bfmu\in\mathcal R\cap\Z^2\\
m\mid\nu_1,\ m^2\mid\nu_2}}
W_r(\bfmu-\bflam)
\ll \frac{r^2}{m^3}+\frac{Nr}{m}+N^2.
\]
It follows from \eqref{eq:Psi-def} and the divisor bound that,
for every $\bflam\in\Z^2$,
\[
\begin{aligned}
\Psi_N(r;2aG_{\bfh},\bflam)
&\ll r^{\eps-1}(r,L(\bfh))
\sum_{m^2\mid r}m
\left(\frac{r^2}{m^3}+\frac{Nr}{m}+N^2\right)\\
&\ll r^\eps(r,L(\bfh))
\bigl(r+N+N^2r^{-1/2}\bigr).
\end{aligned}
\]

For each nonzero integer $n$ with $|n|\ll_F H$, the equation
$L(\bfh)=n$ has $O_F(H)$ solutions with $|\bfh|<H$.
Thus, by the divisor bound,
\[
\begin{aligned}
\sum_{\substack{|\bfh|<H\\L(\bfh)\ne0}}(r,L(\bfh))
&\ll_F H\sum_{1\le n\ll_F H}(r,n)\\
&\ll_F H\sum_{d\mid r}d
\sum_{\substack{1\le n\ll_F H\\d\mid n}}1
\ll H^2r^\eps.
\end{aligned}
\]
Applying \cref{lem:B2} with
$\bflam=2aq_1^2\nabla F(\bfh)$ now gives the first term in
the required bound.

Finally put $\bfv=(-\ell_2,\ell_1)$.
The shifts with $L(\bfh)=0$ are
$\bfh=t\bfv$.  Here $M(\bfv)=D$ and $\nabla F(t\bfv)=0$, while
$G_{\bfv}\circ V=cD^4X^3$.
The calculation in the proof of \cref{lem:singular-average}
applies to this fixed cubic: cancelling $d=(r,t)$ in
\eqref{eq:Th-def} and applying \cref{lem:completion}\textup{(ii)}
at modulus $r/d$ gives
\[
\max_{\calB'}|T_{t\bfv}(\calB')|
\ll r^\eps\bigl(Nr^{1/2}+N^2r^{-1/4}(r,t)^{1/4}\bigr).
\]
By the divisor bound,
\[
\sum_{1\le t<H}(r,t)^{1/4}\ll Hr^\eps.
\]
Hence these shifts satisfy the second term in the statement.
\end{proof}

We now combine the shift estimates with \cref{lem:A2}.
The terms coming from the zero shift and from the leading
completion bound are $N^2q_1^2$ and $N^2q_2$.  Since
$q_1q_2\asymp N^2$, these are balanced by
$q_1\asymp N^{2/3}$ and $q_2\asymp N^{4/3}$.

\begin{proof}[Proof of \cref{thm:main}]
If every linear factor of $F$ occurs at most twice, let
$\mathcal P_F$ be the set chosen after
\cref{lem:polar-discriminant}.  Otherwise use the set in
\cref{lem:triple-factor-average}.
Choose the approximation in \cref{thm:smooth} with $Q=N^2$.
Use \cref{cor:factorization} with $\beta=1/3$ and the fixed
set $\mathcal P_F$.  Then $|\delta|N^4\ll1$,
$q_1=N^{2/3+O(\rho)}$, and $q_2=N^{4/3+O(\rho)}$.
In particular $q_1\le N/2$ for all sufficiently large $N$.

Apply \cref{lem:singular-average,lem:nonsingular-average} in
the first case and \cref{lem:triple-factor-average} in the
second.  Since $N\le Nq_2^{1/4}$, both cases give, on using
\cref{lem:A2} and $q_1H\asymp N$,
\begin{align*}
|S_F(\a;N)|^2
&\ll_{F,\a,\eps} N^\eps
\Bigl(N^2q_1^2+N^2q_2+N^3q_2^{1/4}+N^4q_2^{-1/2}\\
&\hspace{42mm}
+\frac{N^2}{H}
\bigl(q_2+Nq_2^{1/2}+N^2q_2^{-1/4}\bigr)\Bigr).
\end{align*}
Substituting $H\asymp N/q_1$ and the chosen sizes of $q_1,q_2$
gives $|S_F(\a;N)|^2\ll_{F,\a,\eps}
N^{10/3+O(\rho)+\eps}$.
Choose $\rho$ and the auxiliary epsilon small enough in terms
of the epsilon in the Theorem, and take square roots.
\end{proof}

\section{Bounds for complete binary cubic Weyl sums}\label{sec:complete}

Throughout this section we assume $p\ge5$.  We now prove
\cref{thm:complete}.
We reduce the problem to counting solutions of the two derivative
congruences. We count these by grouping the pairs $(x,y)$ according
to the values of $(G_x,G_y)$.

We set $\Sigma(1;G,\bfmu)=1$.  For $a\in\Z_p$, we define
$\e_{p^k}(a)$ by reducing $a$ modulo $p^k$.

We first bound a complete Weyl sum by the number of points where the gradient vanishes
modulo $p^k$.

\begin{lemma}\label{lem:stationary}
Let $n\ge1$ be an integer and let $f\in\Z_p[x_1,\ldots,x_n]$ have
degree at most three.  Put
\[
S_k(f)=\sum_{\bfs\bmod p^k}\e_{p^k}(f(\bfs)),\qquad
N_j(f)=\{\bfs\bmod p^j:\nabla f(\bfs)\equiv0\pmod{p^j}\}.
\]
For $m\ge1$, one has
\begin{equation}\label{eq:half-modulus}
|S_{2m}(f)|\le p^{nm}\#N_m(f),\qquad
|S_{2m+1}(f)|\le p^{n(2m+1)/2}
\bigl(\#N_m(f)\,\#N_{m+1}(f)\bigr)^{1/2}.
\end{equation}
\end{lemma}

\begin{proof}
This is standard and follows from the argument in the proof of
\cite[Lemma~4.1]{Vaughan1997}, applied in several variables.

For the even modulus, write $\bfs=\mathbf a+p^m\mathbf z$,
where $\mathbf a,\mathbf z\bmod p^m$.  Taylor expansion and
orthogonality give
\begin{align*}
S_{2m}(f)
&=\sum_{\mathbf a\bmod p^m}\e_{p^{2m}}(f(\mathbf a))
\sum_{\mathbf z\bmod p^m}
\e_{p^m}(\nabla f(\mathbf a)\cdot\mathbf z)\\
&=p^{nm}\sum_{\mathbf a\in N_m(f)}
\e_{p^{2m}}(f(\mathbf a)).
\end{align*}
This proves the first inequality.

For the odd modulus, write
$\bfs=\mathbf a+p^m\mathbf z+p^{m+1}\mathbf u$,
where $\mathbf a,\mathbf u\bmod p^m$ and $\mathbf z\bmod p$.
For $\mathbf a\in N_m(f)$ put
$\mathbf c_{\mathbf a}=p^{-m}\nabla f(\mathbf a)$ and
$B_{\mathbf a}=(\partial_i\partial_j f(\mathbf a))_{i,j}$.
Taylor expansion and orthogonality in the sum over $\mathbf u$ give
\[
S_{2m+1}(f)
=p^{nm}\sum_{\mathbf a\in N_m(f)}
\e_{p^{2m+1}}(f(\mathbf a))Q_{\mathbf a},
\]
where
\[
Q_{\mathbf a}=
\sum_{\mathbf z\bmod p}
\e_p\bigl(\tfrac12\mathbf z^{\mathsf T}B_{\mathbf a}\mathbf z
+\mathbf c_{\mathbf a}\cdot\mathbf z\bigr).
\]
The cubic terms vanish since $3m\ge2m+1$.
Let $L_{\mathbf a}$ count the solutions of
$B_{\mathbf a}\mathbf z+\mathbf c_{\mathbf a}\equiv0\pmod p$,
which are exactly the lifts of $\mathbf a$ satisfying the
derivative congruences modulo $p^{m+1}$.  Orthogonality gives
\[
|Q_{\mathbf a}|^2
=p^n\sum_{\substack{\mathbf h\bmod p\\B_{\mathbf a}\mathbf h=0}}
\e_p(\mathbf c_{\mathbf a}\cdot\mathbf h)
=p^nL_{\mathbf a},
\]
since $B_{\mathbf a}$ is symmetric.  Hence Cauchy's inequality gives
\[
|S_{2m+1}(f)|
\le p^{nm+n/2}\sum_{\mathbf a\in N_m(f)}\sqrt{L_{\mathbf a}}
\le p^{n(2m+1)/2}\sqrt{\#N_m(f)\,\#N_{m+1}(f)}.
\]
\end{proof}

To apply \cref{lem:stationary} with $f=G+\bfmu\cdot\bfs$,
we must estimate $\#N_j(f)$.

Write $G(x,y)=ax^3+bx^2y+cxy^2+dy^3$ and $\Delta=\disc(G)$.
When $\Delta\ne0$, put
\begin{equation}\label{eq:HG}
H_G(x,y)=(b^2-3ac)x^2+(bc-9ad)xy+(c^2-3bd)y^2.
\end{equation}
Direct substitution gives
\begin{equation}\label{eq:covariants}
G_{xx}G_{yy}-G_{xy}^2=-4H_G,\qquad
\disc(Q_G)=\disc(H_G)=-3\Delta,\qquad
Q_G(G_x,G_y)=H_G^2.
\end{equation}
Here $Q_G=c_GD_G$ is the quadratic form defined in
\eqref{eq:DG}.  In particular,
\[
\disc(D_G)=-3\Delta/c_G^2\ne0.
\]
The last identity will allow us to bound the derivative count
by an average involving the quadratic form $D_G$.

We shall use the elementary bound
\begin{equation}\label{eq:square-root-count}
\#\{x\bmod p^j:Cx^2+v\equiv0\pmod{p^j}\}
\le2\min\{p^{j/2},p^{v_p(v)/2}\}\qquad(p\nmid C).
\end{equation}
If $p^j\mid v$, the number of roots is $p^{\lfloor j/2\rfloor}$.
Otherwise a solution requires $v_p(v)=2s<j$.  Writing $x=p^su$,
Hensel's Lemma gives at most two roots modulo $p^{j-2s}$,
each with $p^s$ lifts modulo $p^{j-s}$.

\begin{lemma}\label{lem:derivative-count}
Let $G$ be an integral binary cubic with $\Delta\ne0$, and let
$\bfmu\in\Z_p^2$.  For $j\ge1$, put
\[
N_j=\{\bfs\bmod p^j:
G_x(\bfs)+\mu_1\equiv G_y(\bfs)+\mu_2\equiv0\pmod{p^j}\}.
\]
Then
\begin{equation}\label{eq:derivative-count}
\#N_j\ll p^{v_p(\Delta)/4}
\min\{p^j,p^{(j+v_p(\bfmu))/2},p^{v_p(D_G(\bfmu))/2}\}.
\end{equation}
\end{lemma}

\begin{proof}
We first consider $\bfs\in N_j$ for which
$p\nmid H_G(\bfs)$ and show that there are at most four.

For any $\mathbf t\in\Z_p^2$, the polynomials $G_x-t_1$ and
$G_y-t_2$ have no common factor.  Indeed, the highest-degree part
of such a factor would divide both $G_x$ and $G_y$, contrary to
$\Delta\ne0$.  The affine B\'ezout theorem
\cite[Theorem~I.2]{Mondal2021} therefore shows that the equations
\[
G_x(\bfr)=t_1,\qquad G_y(\bfr)=t_2
\]
have at most four solutions in $\Z_p^2$.
By Hensel's Lemma, each $\bfs\in N_j$ with
$p\nmid H_G(\bfs)$ lifts uniquely to an exact solution of
\[
G_x(\bfr)=-\mu_1,\qquad G_y(\bfr)=-\mu_2,
\]
and is its reduction modulo $p^j$.  Hence there are at most
four such solutions.

Define
\[
N_j^{\mathrm{sing}}
=\{\bfs\in N_j:p\mid H_G(\bfs)\}.
\]
Since $\#N_j\le4+\#N_j^{\mathrm{sing}}$, it remains to bound
$\#N_j^{\mathrm{sing}}$.  The result is trivial if
$N_j^{\mathrm{sing}}=\varnothing$, so suppose otherwise.

We pass to a larger modulus at which Hensel's Lemma applies to
most of the lifts of these solutions.  The derivative congruences
are still imposed only modulo $p^j$; their values modulo the larger
power of $p$ will be allowed to vary.

Choose an even integer
\[
K\ge2v_p(c_G)+8j+4.
\]
The coefficients of $H_G$ have the same greatest common divisor
as those of $Q_G$.  Thus $p^{-v_p(c_G)}H_G$ has integer
coefficients, with at least one relatively prime to $p$.
After an invertible linear substitution, completing the square
and applying \eqref{eq:square-root-count} gives
\[
\#\{\bfs\bmod p^K:p^{K/2}\mid H_G(\bfs)\}
\le2p^{2K-(K/2-v_p(c_G))/2}
\le\frac2p\,p^{2(K-j)}.
\]
There are $\#N_j^{\mathrm{sing}}\,p^{2(K-j)}$ lifts modulo $p^K$
of the elements of $N_j^{\mathrm{sing}}$.  Each satisfies
\[
G_x(\bfs)+\mu_1\equiv G_y(\bfs)+\mu_2\equiv0\pmod{p^j}.
\]
The preceding bound counts all pairs modulo $p^K$ with
$p^{K/2}\mid H_G(\bfs)$, so it also bounds their number among
these lifts.  Since $p\ge5$ and $N_j^{\mathrm{sing}}\ne\varnothing$,
this number is at most
\[
\frac2p\,p^{2(K-j)}
\le\frac25p^{2(K-j)}
<\frac12\#N_j^{\mathrm{sing}}\,p^{2(K-j)}.
\]
Thus at least half of the lifts satisfy $v_p(H_G(\bfs))<K/2$.
Define
\[
N_{j,K}
=\left\{\bfs\bmod p^K:
\begin{gathered}
G_x(\bfs)+\mu_1\equiv G_y(\bfs)+\mu_2\equiv0\pmod{p^j},\\
1\le v_p(H_G(\bfs))<K/2
\end{gathered}\right\}.
\]
Then
\begin{equation}\label{eq:lift-count-comparison}
\#N_j^{\mathrm{sing}}\le2p^{-2(K-j)}\#N_{j,K}.
\end{equation}

We now group $\bfs\in N_{j,K}$ according to the values
of $G_x(\bfs),G_y(\bfs)$ modulo $p^K$.
Write these values as $\mathbf t=(t_1,t_2)$; necessarily
$\mathbf t\equiv-\bfmu\pmod{p^j}$.
Fix such a $\mathbf t$ and choose a representative in $\Z_p^2$.
For $\bfs\in N_{j,K}$ with these derivative values,
\[
\begin{gathered}
G_x(\bfs)-t_1\equiv G_y(\bfs)-t_2\equiv0\pmod{p^K},\\
H_G(\bfs)^2\not\equiv0\pmod{p^K}.
\end{gathered}
\]
The second condition follows from the restriction
$v_p(H_G(\bfs))<K/2$ in the definition of $N_{j,K}$.
The matrix of first derivatives of $G_x-t_1,G_y-t_2$ satisfies
\[
\det\begin{pmatrix}
G_{xx}(\bfs)&G_{xy}(\bfs)\\
G_{xy}(\bfs)&G_{yy}(\bfs)
\end{pmatrix}
=-4H_G(\bfs).
\]
Since $p\ge5$, the square of this determinant is also nonzero
modulo $p^K$.
The multivariable version of Hensel's Lemma
\cite[(5.21), p.~64]{Greenberg1969}, applied to $G_x-t_1,G_y-t_2$,
therefore gives an exact solution $\bfr\in\Z_p^2$ of
\[
G_x(\bfr)=t_1,\qquad G_y(\bfr)=t_2,
\]
with $v_p(\bfr-\bfs)>v_p(H_G(\bfs))$.
Since $G_x,G_y$ are homogeneous quadratics, their Taylor
expansions give
\[
\begin{pmatrix}
t_1-G_x(\bfs)\\
t_2-G_y(\bfs)
\end{pmatrix}
=
\begin{pmatrix}
G_{xx}(\bfs)&G_{xy}(\bfs)\\
G_{xy}(\bfs)&G_{yy}(\bfs)
\end{pmatrix}
\begin{pmatrix}
r_1-s_1\\
r_2-s_2
\end{pmatrix}
+
\begin{pmatrix}
G_x(\bfr-\bfs)\\
G_y(\bfr-\bfs)
\end{pmatrix}.
\]
Both entries on the left have valuation at least $K$, while
each quadratic remainder has valuation at least
$2v_p(\bfr-\bfs)$.  The inverse matrix has integral numerators
and denominator $-4H_G(\bfs)$.  Multiplying this vector identity
by the inverse matrix and isolating $\bfr-\bfs$ therefore gives
\[
v_p(\bfr-\bfs)\ge
\min\{K-v_p(H_G(\bfs)),
2v_p(\bfr-\bfs)-v_p(H_G(\bfs))\} = K-v_p(H_G(\bfs)).
\]
Consequently,
\[
\bfr\equiv\bfs\pmod{p^{K-v_p(H_G(\bfs))}},
\qquad v_p(H_G(\bfr))=v_p(H_G(\bfs)).
\]

For each exact solution $\bfr$ obtained above, we now bound
\[
\#\left\{\bfs\in N_{j,K}:
\begin{gathered}
G_x(\bfs)\equiv t_1,\quad G_y(\bfs)\equiv t_2\pmod{p^K},\\
\bfs\equiv\bfr\pmod{p^{K-v_p(H_G(\bfr))}}
\end{gathered}\right\}.
\]
Write
\[
\bfs=\bfr+p^{K-v_p(H_G(\bfr))}\mathbf u,
\qquad \mathbf u=(u_1,u_2)\bmod p^{v_p(H_G(\bfr))}.
\]
By Taylor expansion, the congruences
$G_x(\bfs)\equiv t_1$, $G_y(\bfs)\equiv t_2\pmod{p^K}$
are equivalent to
\[
\begin{aligned}
G_{xx}(\bfr)u_1+G_{xy}(\bfr)u_2
&\equiv0\pmod{p^{v_p(H_G(\bfr))}},\\
G_{xy}(\bfr)u_1+G_{yy}(\bfr)u_2
&\equiv0\pmod{p^{v_p(H_G(\bfr))}}.
\end{aligned}
\]
The determinant of the coefficient matrix has valuation
$v_p(H_G(\bfr))$.  Invertible row and column operations over
$\Z_p$ reduce these congruences to
\[
p^{a_1}z_1\equiv0,\qquad p^{a_2}z_2\equiv0
\pmod{p^{v_p(H_G(\bfr))}},
\qquad a_1+a_2=v_p(H_G(\bfr)).
\]
There are $p^{a_1}$ choices for $z_1$ and $p^{a_2}$ choices
for $z_2$ modulo $p^{v_p(H_G(\bfr))}$, so the displayed set
has at most $p^{v_p(H_G(\bfr))}$ elements.

To sum over $\mathbf t$, we now use \eqref{eq:covariants}, which gives
\[
v_p(Q_G(\mathbf t))=2v_p(H_G(\bfr))<K.
\]
There are at most four exact solutions for each $\mathbf t$.
If $p^K\mid Q_G(\mathbf t)$, the preceding identity shows that
no $\bfs\in N_{j,K}$ has these derivative values.  Thus
\[
\begin{aligned}
&\#\{\bfs\in N_{j,K}:
G_x(\bfs)\equiv t_1,\ G_y(\bfs)\equiv t_2\pmod{p^K}\}\\
&\qquad\le
\begin{cases}
4p^{v_p(Q_G(\mathbf t))/2},&v_p(Q_G(\mathbf t))<K,\\
0,&v_p(Q_G(\mathbf t))\ge K
\end{cases}
\ \le4p^{\min(K,v_p(Q_G(\mathbf t)))/2}.
\end{aligned}
\]
The expression on the right depends only on $\mathbf t\bmod p^K$
and is finite even when $Q_G(\mathbf t)=0$.
Summing over the possible derivative values therefore gives
\begin{equation}\label{eq:derivative-average-reduction}
\begin{aligned}
\#N_{j,K}
&\le4\sum_{\substack{\mathbf t\bmod p^K\\
\mathbf t\equiv-\bfmu\bmod p^j}}
p^{\min(K,v_p(Q_G(\mathbf t)))/2}\\
&\le4p^{v_p(c_G)/2}
\sum_{\substack{\mathbf t\bmod p^K\\
\mathbf t\equiv-\bfmu\bmod p^j}}
p^{\min(K,v_p(D_G(\mathbf t)))/2}.
\end{aligned}
\end{equation}

By \eqref{eq:lift-count-comparison} and
\eqref{eq:derivative-average-reduction}, the required bound for
$\#N_j$ will follow if we show that
\begin{equation}\label{eq:quadratic-sum-target}
\begin{aligned}
&\sum_{\substack{\mathbf t\bmod p^K\\
\mathbf t\equiv-\bfmu\bmod p^j}}
p^{\min(K,v_p(D_G(\mathbf t)))/2}\\
&\qquad\ll p^{2(K-j)+v_p(\disc D_G)/4}
\min\{p^j,p^{(j+v_p(\bfmu))/2},p^{v_p(D_G(\bfmu))/2}\}.
\end{aligned}
\end{equation}
Here we use $\disc(D_G)=-3\Delta/c_G^2$.
Changing $\mathbf t$ to $-\mathbf t$ replaces the residue
condition by $\mathbf t\equiv\bfmu\pmod{p^j}$, since
$D_G(-\mathbf t)=D_G(\mathbf t)$.

We put the quadratic form $D_G$ into a form for which this
average can be calculated directly.  Interchanging its variables
or replacing the second variable by the sum of the two variables
makes the coefficient of $X^2$ relatively prime to $p$.
Completing the square then gives
$A(X^2-\beta Y^2)$, where $p\nmid A$ and
$v_p(\beta)=v_p(\disc D_G)$.
These changes permute the summation vectors modulo $p^K$.
Write $\boldsymbol\nu=(\nu_1,\nu_2)$ for the image of $\bfmu$.
Then
\[
v_p(\boldsymbol\nu)=v_p(\bfmu),\qquad
v_p(\nu_1^2-\beta\nu_2^2)=v_p(D_G(\bfmu)).
\]
Since $p\nmid A$, dividing \eqref{eq:quadratic-sum-target}
by $p^{2(K-j)}$ gives the equivalent bound
\begin{equation}\label{eq:quadratic-average}
\begin{aligned}
&p^{-2(K-j)}
\sum_{\substack{X,Y\bmod p^K\\
(X,Y)\equiv\boldsymbol\nu\pmod{p^j}}}
p^{\min(v_p(X^2-\beta Y^2),K)/2}\\
&\qquad\ll p^{v_p(\beta)/4}
\min\{p^j,p^{(j+v_p(\boldsymbol\nu))/2},
p^{v_p(\nu_1^2-\beta\nu_2^2)/2}\}.
\end{aligned}
\end{equation}
We prove this for every $K\ge2j$.

We begin with two averages over all residues that will be used
after taking out common powers of $p$.
For every $R\ge0$, we claim that
\begin{equation}\label{eq:quadratic-full-average}
\begin{aligned}
p^{-R}\sum_{z\bmod p^R}
p^{\min(v_p(z^2-\beta),R)/2}
&\ll1+v_p(\beta),\\
p^{-2R}\sum_{X,Y\bmod p^R}
p^{\min(v_p(X^2-\beta Y^2),R)/2}
&\ll1+v_p(\beta)\ll p^{v_p(\beta)/4}.
\end{aligned}
\end{equation}
The case $R=0$ is trivial, so suppose $R\ge1$.
For the first bound, \eqref{eq:square-root-count} gives
\[
\begin{aligned}
p^{-R}\sum_{z\bmod p^R}p^{\min(v_p(z^2-\beta),R)/2}
&\le1+\sum_{r=1}^Rp^{-r/2}
\#\{z\bmod p^r:z^2\equiv\beta\pmod{p^r}\}\\
&\ll1+\min\{R,v_p(\beta)\}
 +\sum_{v_p(\beta)<r\le R}p^{(v_p(\beta)-r)/2}\\
&\ll1+v_p(\beta).
\end{aligned}
\]

For the second bound, splitting according to whether $p\mid Y$
and substituting $X=zY$ when $p\nmid Y$ gives
\[
\begin{aligned}
\sum_{\substack{X,Y\bmod p^R\\p\nmid(X,Y)}}
p^{\min(v_p(X^2-\beta Y^2),R)/2}
&=\sum_{\substack{X,Y\bmod p^R\\p\nmid X,\ p\mid Y}}1\\
&\quad+\sum_{\substack{Y\bmod p^R\\p\nmid Y}}
\sum_{z\bmod p^R}p^{\min(v_p(z^2-\beta),R)/2}\\
&\ll p^{2R}+p^R\cdot p^R(1+v_p(\beta))\\
&\ll p^{2R}(1+v_p(\beta)).
\end{aligned}
\]

We now separate all pairs according to $m=v_p(X,Y)$.
For $2m<R$, write $X=p^mx$ and $Y=p^my$, where
$x,y\bmod p^{R-m}$ and $p\nmid(x,y)$.  The preceding bound gives
\[
\begin{aligned}
p^{-2R}\sum_{\substack{X,Y\bmod p^R\\v_p(X,Y)=m}}
p^{\min(v_p(X^2-\beta Y^2),R)/2}
&\ll p^{-2R}p^mp^{2(R-m)}(1+v_p(\beta))\\
&=p^{-m}(1+v_p(\beta)).
\end{aligned}
\]
The remaining pairs have both coordinates divisible by
$p^{\lceil R/2\rceil}$, so
\[
v_p(X^2-\beta Y^2)\ge2\lceil R/2\rceil\ge R,
\qquad
p^{\min(v_p(X^2-\beta Y^2),R)/2}=p^{R/2}.
\]
There are $p^{2(R-\lceil R/2\rceil)}$ such pairs.  Consequently,
\[
\begin{aligned}
p^{-2R}\sum_{X,Y\bmod p^R}p^{\min(v_p(X^2-\beta Y^2),R)/2}
&\ll(1+v_p(\beta))\sum_{0\le m<R/2}p^{-m}
 +p^{R/2-2\lceil R/2\rceil}\\
&\ll1+v_p(\beta)\ll p^{v_p(\beta)/4}.
\end{aligned}
\]

We first prove \eqref{eq:quadratic-average} when
$p\nmid\boldsymbol\nu$.  If $p\mid\nu_2$, then $p\nmid\nu_1$
and every summand in that average is one, which gives the result.
Otherwise $p\nmid\nu_2$.  Put $z_0=\nu_1/\nu_2$ and substitute
$X=zY$.  The average on the left of \eqref{eq:quadratic-average}
becomes
\[
J=p^{-(K-j)}
\sum_{\substack{z\bmod p^K\\z\equiv z_0\pmod{p^j}}}
p^{\min(v_p(z^2-\beta),K)/2}.
\]
Since $v_p(\nu_1^2-\beta\nu_2^2)=v_p(z_0^2-\beta)$, it remains
to prove
\begin{equation}\label{eq:quadratic-local-average}
J\ll p^{v_p(\beta)/4}
\min\{p^{j/2},p^{v_p(z_0^2-\beta)/2}\}.
\end{equation}

Suppose first that $v_p(z_0)<j$.  Squaring converts the restricted
average into one involving a linear polynomial.  Indeed, for
$z,z'\equiv z_0\pmod{p^j}$,
\[
v_p(z^2-z'^2)=v_p(z-z')+v_p(z_0).
\]
Thus each value of $z^2\bmod p^K$ has $p^{v_p(z_0)}$ preimages,
so there are $p^{K-j-v_p(z_0)}$ such values.  They are therefore
precisely the residues congruent to $z_0^2$ modulo
$p^{j+v_p(z_0)}$.  Hence
\[
J=p^{-(K-j-v_p(z_0))}
\sum_{u\bmod p^{K-j-v_p(z_0)}}
p^{\min(v_p(z_0^2-\beta+p^{j+v_p(z_0)}u),K)/2}.
\]
If $v_p(z_0^2-\beta)<j+v_p(z_0)$, every summand is
$p^{v_p(z_0^2-\beta)/2}$.  Otherwise translate $u$ by
$(z_0^2-\beta)/p^{j+v_p(z_0)}$ to obtain
\[
\begin{aligned}
J
&=p^{(j+v_p(z_0))/2}p^{-(K-j-v_p(z_0))}
\sum_{v\bmod p^{K-j-v_p(z_0)}}
p^{\min(v_p(v),K-j-v_p(z_0))/2}\\
&\le p^{(j+v_p(z_0))/2}
\sum_{r=0}^{K-j-v_p(z_0)}p^{-r/2}
\ll p^{(j+v_p(z_0))/2}.
\end{aligned}
\]
Combining the two cases gives
\[
J\ll p^{\min\{v_p(z_0^2-\beta),\,j+v_p(z_0)\}/2}.
\]
This gives the term involving $v_p(z_0^2-\beta)$ in
\eqref{eq:quadratic-local-average}.
If $2v_p(z_0)\le v_p(\beta)$, it also gives
$J\ll p^{j/2+v_p(\beta)/4}$.
If $2v_p(z_0)>v_p(\beta)$, then
$v_p(z_0^2-\beta)=v_p(\beta)<2j$, which gives the same bound.
This proves \eqref{eq:quadratic-local-average} when $v_p(z_0)<j$.

If $v_p(z_0)\ge j$, the summation condition is simply
$p^j\mid z$.  When $v_p(\beta)<2j$, every summand is
$p^{v_p(\beta)/2}$ and $v_p(z_0^2-\beta)=v_p(\beta)$,
which gives \eqref{eq:quadratic-local-average}.
When $v_p(\beta)\ge2j$, substitute $z=p^ju$ and apply the first
estimate in \eqref{eq:quadratic-full-average} to $\beta/p^{2j}$:
\[
\begin{aligned}
J
&=p^jp^{-(K-2j)}
\sum_{u\bmod p^{K-2j}}
p^{\min(v_p(u^2-\beta/p^{2j}),K-2j)/2}\\
&\ll p^j(1+v_p(\beta)-2j)
\ll p^{j/2+v_p(\beta)/4}.
\end{aligned}
\]
Here $v_p(z_0^2-\beta)\ge2j$, so this proves
\eqref{eq:quadratic-local-average} and completes
\eqref{eq:quadratic-average} when $p\nmid\boldsymbol\nu$.

We now remove the assumption $p\nmid\boldsymbol\nu$.
If $\boldsymbol\nu\equiv0\pmod{p^j}$, put $X=p^jx$ and
$Y=p^jy$.  The second estimate in \eqref{eq:quadratic-full-average}
gives
\[
\begin{aligned}
&p^{-2(K-j)}
\sum_{\substack{X,Y\bmod p^K\\p^j\mid X,\ p^j\mid Y}}
p^{\min(v_p(X^2-\beta Y^2),K)/2}\\
&\qquad=p^jp^{-2(K-2j)}
\sum_{x,y\bmod p^{K-2j}}
p^{\min(v_p(x^2-\beta y^2),K-2j)/2}
\ll p^{j+v_p(\beta)/4}.
\end{aligned}
\]
Here the minimum in \eqref{eq:quadratic-average} is $p^j$.

Otherwise $0<v_p(\boldsymbol\nu)<j$, so write
\[
\boldsymbol\nu=p^{v_p(\boldsymbol\nu)}\boldsymbol\xi,
\qquad (X,Y)=p^{v_p(\boldsymbol\nu)}(x,y).
\]
Then $p\nmid\boldsymbol\xi$.  Homogeneity gives a factor
$p^{v_p(\boldsymbol\nu)}$ in the normalized average, and leaves
the modulus and congruence condition
\[
x,y\bmod p^{K-2v_p(\boldsymbol\nu)},\qquad
(x,y)\equiv\boldsymbol\xi\pmod{p^{j-v_p(\boldsymbol\nu)}}.
\]
Applying the case just proved, with $K,j$ replaced by
$K-2v_p(\boldsymbol\nu),j-v_p(\boldsymbol\nu)$, therefore gives
\[
\begin{aligned}
&p^{-2(K-j)}
\sum_{\substack{X,Y\bmod p^K\\
(X,Y)\equiv\boldsymbol\nu\pmod{p^j}}}
p^{\min(v_p(X^2-\beta Y^2),K)/2}\\
&\qquad\ll p^{v_p(\boldsymbol\nu)+v_p(\beta)/4}
\min\{p^{(j-v_p(\boldsymbol\nu))/2},
p^{v_p(\xi_1^2-\beta\xi_2^2)/2}\}\\
&\qquad=p^{v_p(\beta)/4}
\min\{p^{(j+v_p(\boldsymbol\nu))/2},
p^{v_p(\nu_1^2-\beta\nu_2^2)/2}\}.
\end{aligned}
\]
Since $v_p(\boldsymbol\nu)<j$, this proves
\eqref{eq:quadratic-average}, and hence
\eqref{eq:quadratic-sum-target}.
Inserting this estimate into \eqref{eq:derivative-average-reduction}
and using \eqref{eq:lift-count-comparison}, we obtain
\[
\begin{aligned}
\#N_j
&\le4+\#N_j^{\mathrm{sing}}
\le4+2p^{-2(K-j)}\#N_{j,K}\\
&\ll p^{v_p(\Delta)/4}
\min\{p^j,p^{(j+v_p(\bfmu))/2},p^{v_p(D_G(\bfmu))/2}\}.
\end{aligned}
\]
This proves \eqref{eq:derivative-count}.
\end{proof}

\begin{proof}[Proof of \cref{thm:complete}\textup{(i)}]
We first work modulo $p^k$.  Apply \cref{lem:stationary} to
$f=G+\bfmu\cdot\bfs$, and insert \eqref{eq:derivative-count}.
For both parities of $k\ge2$ this gives
\[
|\Sigma(p^k;G,\bfmu)|
\ll p^{k+v_p(\Delta)/4}
\min\{p^{k/2},p^{k/4+v_p(\bfmu)/2},p^{v_p(D_G(\bfmu))/2}\}.
\]
The two valuations in the minimum may be capped at $k$.
Using the trivial bound $p^{2k}$ when $v_p(\Delta)\ge4k$, we obtain
\begin{equation}\label{eq:local-bound}
|\Sigma(p^k;G,\bfmu)|
\ll p^{k+\min(k,v_p(\Delta)/4)}
\min\{p^{k/4}(p^k,\bfmu)^{1/2},
(p^k,D_G(\bfmu))^{1/2}\}.
\end{equation}

For $k=1$ and $p\nmid\Delta$, Deligne's bound
\cite{Deligne1974} gives $4p$.
If $G$ is nonzero modulo $p$ and has a repeated factor modulo $p$,
a linear substitution invertible modulo $p$ puts it in the form
$AX^2Y$ or $AX^3$ modulo $p$, with $p\nmid A$.
In the first case, writing the new
linear coefficients as $(u,v)$ gives
\[
\sum_{X,Y\bmod p}\e_p(AX^2Y+uX+vY)
=p\sum_{\substack{X\bmod p\\AX^2+v=0}}\e_p(uX) \ll p.
\]
In the cube case, summing over $Y$
and applying the one-variable Weil bound gives $2p^{3/2}$.
Lifting the substitution to $\Z_p$ and using the cubic
discriminant formula gives $v_p(\Delta)\ge2$ in this case.
If every coefficient of $G$ is divisible by $p$, then $v_p(\Delta)\ge4$
and the trivial bound is $p^2$.  Thus \eqref{eq:local-bound}
also holds for $k=1$.

The Chinese remainder theorem and the divisor bound now give
\eqref{eq:complete-bound}.
\end{proof}

\begin{proof}[Proof of \cref{thm:complete}\textup{(ii)}]
The repeated linear factor of $G$ is rational, as is the
remaining simple factor when there is one.  A rational linear
substitution therefore puts $G$ in the form $aX^2Y$ or $aX^3$.
Multiplying its matrix by a common denominator gives an integer
matrix $V$ with $\det V\ne0$ and
\[
G\circ V=AX^2Y\quad\hbox{or}\quad G\circ V=AX^3,
\qquad A\in\Z\setminus\{0\}.
\]
For $(r,6A\det V)=1$, the substitution
$\bfs=V(X,Y)^{\mathsf T}$ permutes the pairs modulo $r$,
and the linear term becomes $MX+LY$, where
$M=M_G(\bfmu)$ and $L=L_G(\bfmu)$.

In the first case, summing over $Y$ gives
\[
\Sigma(r;G,\bfmu)
=r\sum_{\substack{X\bmod r\\AX^2+L\equiv0\pmod r}}
\e_r(MX).
\]
By \eqref{eq:square-root-count}, the Chinese remainder theorem
and the divisor bound, the number of terms is
$O(r^\eps(r,L)^{1/2})$.
This proves \eqref{eq:complete-double}.
In the second case the sum factors as
\[
\Sigma(r;G,\bfmu)
=\left(\sum_{X\bmod r}\e_r(AX^3+MX)\right)
\left(\sum_{Y\bmod r}\e_r(LY)\right).
\]
The second factor is $r\mathbf1_{r\mid L}$.
Applying \cite[Lemma~4.1 and Theorem~7.1]{Vaughan1997} to the
first proves \eqref{eq:complete-cube}.
\end{proof}

The following estimate retains the dependence on the coefficient
of $X^2Y$.

\begin{lemma}\label{lem:mixed-cubic}
Let $r\ge1$ satisfy $(r,6)=1$, let $A,B\in\Z$, and let
$\bfmu=(\mu_1,\mu_2)\in\Z^2$.  Then
\begin{equation}\label{eq:mixed-cubic-bound}
|\Sigma(r;AX^3+BX^2Y,\bfmu)|
\ll r^{1+\eps}(r,B)(r,\mu_1^2,\mu_2)^{1/2}.
\end{equation}
\end{lemma}

\begin{proof}
Summing over $Y$ gives
\[
\Sigma(r;AX^3+BX^2Y,\bfmu)
=r\sum_{\substack{X\bmod r\\BX^2+\mu_2\equiv0\pmod r}}
\e_r(AX^3+\mu_1X).
\]
The case $r=1$ is immediate.  Suppose that $r=p^k$ and put
$t=\min(v_p(B),k)$.  First restrict to $v_p(X)=s<k/2$.
Translation by $p^{k-s}z$ preserves this valuation and leaves
$BX^2$ and $AX^3$ unchanged modulo $p^k$.  Moreover,
\[
\sum_{z\bmod p^s}
\e_{p^k}\bigl(A(X+p^{k-s}z)^3+\mu_1(X+p^{k-s}z)\bigr)
=p^s\e_{p^k}(AX^3+\mu_1X)\mathbf1_{p^s\mid\mu_1}.
\]
Thus the contribution from these $X$ vanishes unless
$p^s\mid\mu_1$.  The congruence also requires
$p^{2s}\mid\mu_2$.

There are at most $2p^{t+s}$ such $X$.  Indeed, if $t+2s<k$,
a solution requires $p^{t+2s}\mid\mu_2$.  Writing $X=p^sZ$
and dividing the congruence by $p^{t+2s}$ gives
\[
p^{-t}B\,Z^2\equiv-\mu_2/p^{t+2s}
\pmod{p^{k-t-2s}}.
\]
There are at most two solutions with $p\nmid Z$ modulo
$p^{k-t-2s}$, each with $p^{t+s}$ lifts modulo $p^{k-s}$.
If $t+2s\ge k$, the trivial count $p^{k-s}\le p^{t+s}$
suffices.

For the remaining $X$, put $\ell=\floor{k/2}$.  Their
contribution to the sum over $X$ is
\[
\sum_{\substack{X\bmod p^k\\p^{k-\ell}\mid X\\
BX^2+\mu_2\equiv0\pmod{p^k}}}
\e_{p^k}(AX^3+\mu_1X)
=p^\ell\mathbf1_{p^\ell\mid\mu_1}
\mathbf1_{p^k\mid\mu_2}.
\]
Combining these estimates gives
\[
|\Sigma(p^k;AX^3+BX^2Y,\bfmu)|
\ll p^{k+t}
\sum_{\substack{0\le s\le\floor{k/2}\\
p^s\mid\mu_1,\ p^{2s}\mid\mu_2}}p^s.
\]
The Chinese remainder theorem and the divisor bound give
\[
\begin{aligned}
|\Sigma(r;AX^3+BX^2Y,\bfmu)|
&\ll r^{1+\eps}(r,B)
\sum_{m^2\mid(r,\mu_1^2,\mu_2)}m\\
&\ll r^{1+\eps}(r,B)(r,\mu_1^2,\mu_2)^{1/2}.
\end{aligned}
\]
This proves \eqref{eq:mixed-cubic-bound}.
\end{proof}

% Earlier supplementary remarks are retained in the revision archive.

\printbibliography

\end{document}